\documentclass[11pt]{article}

\usepackage{vien}
\usepackage{mathtools, graphicx, caption, subcaption, bm, epstopdf, url}
\usepackage{algorithm, enumitem}
\usepackage[top=2.54cm,right=3cm,left=3cm,bottom=2.54cm]{geometry}
\usepackage[final]{changes}

\newcommand{\statrv}{S}
\newcommand{\statval}{s}
\newcommand{\stepsize}{\alpha}
\newcommand{\mmt}{\beta}
\newcommand{\subdiff}{\partial}
\newcommand{\opt}{^\star}
\newcommand{\indcfunc}[1]{\mathsf{I}_{\mc{X}}(#1)}
\newcommand{\proj}[1]{{\mathop{\Pi}}_{\mc{X}} \left(#1\right)}

\title{Convergence of a Stochastic Gradient Method with Momentum for Non-Smooth Non-Convex Optimization}

\author{	
	Vien V.~Mai\footnotemark[1]
	\and  
	Mikael Johansson{\thanks{Division of Decision and Control Systems, School of Electrical Engineering and Computer Science, KTH~Royal Institute of Technology, SE-100 44  Stockholm, Sweden. Emails: \tt\small\{maivv, mikaelj\}@kth.se.}}       
}

\begin{document}
\maketitle

\begin{abstract}
Stochastic gradient methods with momentum are widely used in applications and at the core of optimization subroutines in many popular machine learning libraries. However, their sample complexities have not been obtained for problems beyond those that are convex or smooth. This paper establishes the convergence rate of a stochastic subgradient method with a momentum term of Polyak type for a broad class of non-smooth, non-convex, and constrained optimization problems. Our key innovation is the construction of a special Lyapunov function for which the proven complexity can be achieved without any tuning of the momentum parameter. For smooth problems, we extend the known complexity bound to the constrained case and demonstrate how the unconstrained case can be analyzed under weaker assumptions than the state-of-the-art. Numerical results confirm our theoretical developments.
\end{abstract}

\section{Introduction}

We study the stochastic optimization problem
\begin{equation}
\label{eqn:objective}
		\underset{x\in \mc{X}}{\minimize} \,
	 f(x) := \E_{P}[f(x;\statrv)]
	= \int_{\mc{S}} f(x; \statval) dP(\statval),
\end{equation}
where $\statrv\sim P$ is a random variable;  $f(x;\statval)$ is the instantaneous loss parameterized by $x$ on a sample $\statval\in\mc{S}$; and $\mc{X}\subseteq\R^n$ is a closed convex set. 
In this paper, we move beyond convex and/or smooth optimization and consider $f$ that belongs to a broad class of \emph{non-smooth} and \emph{non-convex} functions called $\rho$-weakly convex, meaning that
\begin{align*}
	x \mapsto f(x) + \rho \norm{x}_2^2 \,\,\, \mbox{is convex}.
\end{align*}
This function class  is very rich and important in optimization \cite{Roc81,Via83}. It trivially includes all convex functions, all smooth functions with Lipschitz continuous gradient, and all  additive composite functions of the two former classes. More broadly, it includes all compositions of the form
\begin{align}\label{eq:composite}
	f(x) = h(c(x)),
\end{align}
where $h: \R^m \to \R$ is convex and $L_h$-Lipschitz and $c: \R^n \to \R^m$ is a smooth map with $L_c$-Lipschitz Jacobian. Indeed, the composite function $f = h\circ c$ is then weakly convex with $\rho=L_h L_c$ \cite{DP19}.
Some representative applications in this problem class include nonlinear least squares \cite{Dru17}, robust phase retrieval \cite{DR18a}, Robust PCA \cite{CLMW11}, robust low rank matrix recovery \cite{CCD19}, optimization of the Conditional Value-at-Risk  \cite{RU00},  graph synchronization \cite{Sin11}, and many others. 
   
Stochastic optimization algorithms for solving \eqref{eqn:objective}, based on random samples $S_k$ drawn from $P$, are of fundamental importance in many applied sciences \cite{Bot10, SDR14}.    
Since the introduction of the classical stochastic (sub)gradient descent method (SGD) in \cite{RM51},
several modifications of SGD have been proposed to improve its practical and theoretical performance. A notable example  is the use of a \emph{momentum} term to construct an update  direction \cite{GB72, Pol87, RS83, Tse98, SMDH13, GL16}. The basis form of such a method (when $\mc{X}=\R^n$) reads: 
\begin{subequations}
	\label{eq:SHB:unconstrained} 
  \begin{align}
  	x_{k+1} &= x_k - \stepsize_k z_k\\\label{eq:SHB:unconstrained:z}
     z_{k+1} &=   \mmt_k g_{k+1} + (1-\mmt_k) z_k,
 \end{align}
\end{subequations}  
where $z_k$ is the search direction,  $g_k$ is a stochastic subgradient; $\stepsize_k$ is the stepsize, and $\mmt_k\in(0,1]$ is the momentum parameter. For instance, the scheme \eqref{eq:SHB:unconstrained} reduces to the stochastic heavy ball method (SHB) \cite{Pol87}:
\begin{align*}	
  	x_{k+1} &= x_k - \eta_k g_k + \lambda_k(x_k-x_{k-1}) ,
\end{align*}
with  $\eta_k=\stepsize_k\mmt_{k-1} $ and $\lambda_k=(1-\mmt_{k-1})\stepsize_k/\stepsize_{k-1}$.  Methods of this type enjoy widespread empirical success in large-scale convex and non-convex optimization, especially in training neural networks, where they have been used to produce several state-of-the-art results on important learning tasks, e.g., \cite{KSH12, SMDH13, HZRS16,  HLVW17}.

Sample complexity, namely,  the number of observations $S_0, \ldots, S_K$ required to reach a desired accuracy $\epsilon$, has been the most widely adapted metric for evaluating the performance of stochastic optimization algorithms.  Although sample complexity results for the standard SGD on problems of the form \eqref{eqn:objective} have been obtained for convex and/or smooth problems \cite{NJLS09,GL13}, much less is known in the non-smooth non-convex case \cite{DD19}. The problem is even more prevalent in momentum-based methods as there is virtually no known complexity results for problems beyond those that are convex or smooth.  

\subsection{Related work}
As many applications in  modern machine learning and signal processing  cannot be captured by convex models, (stochastic) algorithms for solving non-convex problems have been studied extensively. Below we  review some of the topics most closely related to our work. \medskip

\noindent\textbf{Stochastic weakly convex minimization}\hspace*{0.5em}
Earlier works on this topic date back to Nurminskii who showed subsequential convergence to stationary points  for the subgradient method applied to {deterministic} problems \cite{Nur73}. The work \cite{Rus87} proposes a stochastic gradient averaging-based method and shows the first almost sure convergence for this problem class.  Basic sufficient conditions for convergence of stochastic projected subgradient methods is established in \cite{EN98}. Thanks to the recent advances in statistical learning and signal processing, the problem class has been reinvigorated with several new theoretical results and practical applications (see, e.g., \cite{DR18a,DR16,DG19,DD19} and references therein). In particular, based on the theory of non-convex differential inclusions, almost sure convergence is derived in \cite{DR18b} for a collection of model-based minimization strategies, albeit no rates of convergence are given. 
An important step toward \emph{non-asymptotic} convergence of stochastic methods is made in \cite{DG19}. There, the authors employ a  proximal point technique for which they can show the sample complexity $O(1/\epsilon^2)$ with a certain stationarity measure.  Later, the work \cite{DD19} shows that the (approximate) proximal point step in \cite{DG19} is not necessary and establishes the similar complexity for a class of model-based methods including  the standard SGD.
We also note that there has been a large volume of works in smooth non-convex optimization, e.g., \cite{GL13, GL16}. \medskip

\noindent\textbf{{Stochastic momentum for non-convex functions}}\hspace*{0.5em}
Optimization algorithms based on momentum averaging techniques go back to Polyak \cite{Pol64} who proposed the heavy ball method. In \cite{Nes83}, Nesterov introduced the accelerated gradient method and showed its optimal iteration complexity for the minimization of smooth convex functions. In the last few decades, research on accelerated first-order methods has exploded both in theory and in practice \cite{Bec17,Nes18,Buc15}. 
The effectiveness of such techniques in the deterministic context has inspired researchers to  incorporate momentum terms into stochastic optimization algorithms \cite{Pol87, Rus87, Tse98, SMDH13, GL16}. 
Despite evident success, especially, in training neural networks \cite{KSH12, SMDH13, HZRS16, ZK16, HLVW17}, the theory for stochastic momentum methods is not as clear as its deterministic counterpart (cf. \cite{GLZX19}).  As a result, there has been a growing interest in obtaining convergence guarantees for those methods under noisy gradients \cite{HPK09,GLZX19,YYL18,GPS18,GL16}. In  non-convex optimization, almost certain convergence of Algorithm~\eqref{eq:SHB:unconstrained} for \emph{smooth and unconstrained} problems is derived in \cite{RS83}. Under the \emph{bounded gradient} hypothesis, the convergence rate of the same algorithm has been established in \cite{YYL18}. The work \cite{GRW18} obtains the complexity of a  gradient averaging-based method for \emph{constrained} problems.  In \cite{GL16}, the authors study a variant of Nesterov acceleration and establish a similar complexity for smooth and unconstrained problems, while for the constrained case, a mini-batch of samples at each iteration is required to guarantee convergence.

\subsection{Contributions}

Minimization of weakly convex functions has been a challenging task, especially for stochastic problems, as the objective is neither smooth nor convex.  With the recent breakthrough in \cite{DD19}, this problem class has been the widest one for which provable sample complexity of the standard SGD is known. It is thus intriguing to ask whether such a result can also be obtained for momentum-based methods. The work in this paper aims to address this question. To that end, we make the following contributions:
\begin{itemize}[leftmargin=0.5cm]
\item We establish the sample complexity of a stochastic subgradient method with momentum of Polyak type for a broad class of non-smooth, non-convex, and constrained optimization problems. Concretely, using a special Lyapunov function, we show the complexity $O(1/\epsilon^2)$ for the minimization of weakly convex functions. The proven complexity is attained in a parameter-free and single time-scale fashion, namely, the stepsize and the momentum constant are independent of any problem parameters and they have the same scale with respect to the iteration count. To the best of our knowledge, this is the first complexity guarantee for a stochastic method with momentum on non-smooth and non-convex problems. 

\item We also study the sample complexity of the considered algorithm for smooth and constrained optimization problems. 
Note that even in this setting, no complexity guarantee of SGD with Polyak momentum has been established before. Under a bounded gradient assumption, we obtain a similar  $O(1/\epsilon^2)$ complexity without the need of forming a batch of samples at each iteration, which is commonly required for \emph{constrained} non-convex stochastic optimization \cite{ GLZ16}. We then demonstrate how the unconstrained case can be analyzed without the above assumption.

\end{itemize}

\noindent Interestingly,  the stated result is achieved in the regime where $\mmt$ can be as small as $O(1/\sqrt{K})$, i.e., one can put much more weight to the momentum term than the fresh subgradient in a search direction. This complements the complexity of SGD attained as $\mmt\to 1$. Note that the worst-case complexity $O(1/\epsilon^2)$ is unimprovable in the smooth and unconstrained case \cite{ACD19}. 

\section{Background}

In this section, we first introduce the notation and then provide the necessary preliminaries for the paper.

For any $x,y\in\R^n$, we denote by $\InP{x}{y}$ the Euclidean inner product of $x$ and $y$. We use $\ltwo{\cdot}$ to denote the Euclidean norm. For a closed and convex set $\mc{X}$,  ${\mathop{\Pi}}_{\mc{X}}$ denotes the orthogonal projection onto $\mc{Z}$, i.e., $y=\proj{x}$ if $y\in\mc{X}$ and $\ltwo{y-x}=\min_{z\in\mc{X}}\ltwo{z-x}$; $\indcfunc{\cdot}$ denotes the indicator function of $\mc{X}$, i.e., $\indcfunc{x}=0$ if $x\in\mc{X}$ and $+\infty$ otherwise. Finally, we denote by $\mc{F}_k := \sigma(\statrv_0, \ldots, \statrv_k)$ the $\sigma$-field generated by the first $k+1$ random variables $\statrv_0, \ldots,\statrv_k$.

For a function $f: \R^n \to \R\cup\{+\infty\}$, the Fr\'{e}chet subdifferential of $f$ at $x$, denoted by $\subdiff f(x)$, consists of all vectors $g\in\R^n$ such that
\begin{align*}
	f(y) \geq f(x) + \InP{g}{y-x} + o(\norm{y-x}) \,\,\, \mbox{as}\,\,\, y\to x.
\end{align*}

The Fr\'{e}chet and conventional subdifferentials coincide for convex functions, while for smooth functions $f$,  $\subdiff f(x)$ reduces to the gradient $\{\grad{f(x)}\}$. 
A point $x\in\R^n$ is said to be \emph{stationary} for problem \eqref{eqn:objective} if $0\in \subdiff{f(x)+\subdiff\indcfunc{x}}$.

The following lemma collects standard properties of weakly convex functions \cite{Via83}.
\begin{lemma}[Weak convexity]
Let $f: \R^n \to \R\cup\{+\infty\}$ be a $\rho$-weakly convex function. Then the following hold:
\begin{description}[leftmargin=0.45cm]
\item[1.] For any $x, y\in \R^n$ with $g\in \subdiff f(x)$, we have
\begin{align*}
	f(y) \geq f(x) + \InP{g}{y-x} - \frac{\rho}{2} \norm{y-x}_2^2. 
\end{align*}
\item[2.] For all $x, y\in \R^n$, $\alpha\in [0,1]$, and $z=\alpha x + (1-\alpha)y$: 
\begin{align*}
	f(z) 
	&\leq 
		\alpha f(x) + (1-\alpha)f(y) 
		+ \frac{\rho \alpha(1-\alpha)}{2}\norm{x-y}_2^2.
\end{align*}
\end{description}
\end{lemma}
Weakly convex functions  admit an implicit smooth approximation through the Moreau envelope:
\begin{align}\label{eq:moreau:env}
	f_\lambda(x) = \inf_{y\in\R^n} \left\{f(y)+ \frac{1}{2\lambda} \ltwo{x-y}^2\right\}.
\end{align}
For small enough $\lambda$, the point achieving $f_\lambda(x)$ in \eqref{eq:moreau:env}, denoted by $\prox{\lambda f}{x} $,  is unique and given by:
\begin{align}
	\label{eq:prox}
	\prox{\lambda f}{x} = \argmin_{y\in\R^n} \left\{f(y)+ \frac{1}{2\lambda} \ltwo{x-y}^2\right\}.
\end{align}
The lemma below summarizes two well-known properties of the Moreau envelope and its associated proximal map \cite{HL93}.
\begin{lemma}[Moreau envelope]\label{lem:moreau:env}
Suppose that $f: \R^n \to \R\cup\{+\infty\}$ is a $\rho$-weakly convex function. Then, for a fixed parameter $\lambda^{-1} \geq 2\rho$, the following hold:
\begin{description}[leftmargin=0.45cm]
\item[1.]  $f_\lambda$ is $\mc{C}^1$-smooth with the gradient given by
\begin{align*}
	\grad{f_{\lambda}(x)} = \lambda^{-1}\left(x-\prox{\lambda f}{x} \right).
\end{align*}
\item[2.]  $f_\lambda$ is $\lambda^{-1}$-smooth, i.e., for all $x, y\in \R^n$:
\begin{align*}
			\big| f_\lambda(y) - f_\lambda(x) - \InP{\grad{f_{\lambda}(x)}}{y-x} \big|
	\leq 
		\frac{1}{2\lambda}\ltwo{x-y}^2.
\end{align*}
\end{description}
\end{lemma}

\noindent\textbf{Failure of stationarity test}\hspace*{0.5em} A major source of difficulty in convergence analysis of non-smooth optimization methods is the lack of a controllable stationarity  measure.  For smooth functions, it is natural to use the norm of the gradient as a surrogate for near stationarity. However, this rule does not make sense in the non-smooth case, even if the function is convex and its gradient existed at all iterates. For example, the \emph{convex} function $f(x)=|x|$ has $\ltwo{\grad{f(x)}}=1$ at each $x\neq0$, no mater how close $x$ is to the stationary point $x=0$.

To circumvent this difficulty, we  adopt the techniques pioneered in \cite{DD19} for convergence of stochastic methods on weakly convex problems. More concretely, we rely on the connection of the Moreau envelope to (near) stationarity: 
For any $x\in\R^n$, the point $\hat{x}=\prox{\lambda F}{x}$, where $F(x) = f(x) + \indcfunc{x}$, satisfies: 
\begin{align}\label{eq:stationary:measure}
\begin{cases}
	\ltwo{x-\hat{x}} =  \lambda \ltwo{\grad{F_{\lambda}(x)}},\\
	\dist(0, \subdiff F(\hat{x})) \leq \ltwo{\grad{F_{\lambda}(x)}}.
\end{cases}
\end{align}
Therefore, a small gradient $\ltwo{\grad{F_{\lambda}(x)}}$ implies that $x$ is close to a point $\hat{x}\in\mc{X}$ that is near-stationary for $F$. Note that $\hat{x}$ is just a \emph{virtual}  point, there is no need to compute it.

\section{Algorithm and convergence analysis}
We assume that the only access to $f$ is through a stochastic subgradient oracle. In particular, we study algorithms that attempt to solve problem \eqref{eqn:objective} using i.i.d. samples $\statrv_0, \statrv_1, \ldots,\statrv_K\simiid P$. 
\begin{assumption}[Stochastic oracle]\label{assumption:A}
Fix a probability space $(\mc{S}, \mc{F}, P)$. Let $\statrv$ be a sample drawn from $P$ and $f'(x,\statrv) \in \subdiff f(x, \statrv)$. We make the following assumptions:
\begin{description}[leftmargin=0.45cm]
\item[(A1)] For each $x \in \dom(f)$, we have 
\begin{align*}
	\E_{P}\left[f'(x, \statrv)\right] \in \subdiff f(x).
\end{align*}

\item [(A2)] There exists a real $L>0$ such that  for all $x\in \mc{X}$:
\begin{align*}
	\E_{P} \left[\ltwo{f'(x, \statrv)}^2\right] \leq L^2.
\end{align*}
\end{description}
\end{assumption}
The above assumptions are standard in stochastic optimization of \emph{non-smooth} functions (see, e.g., \cite{NJLS09,DD19}).

\noindent\textbf{Algorithm}\hspace*{0.5em} 
To solve problem \eqref{eqn:objective}, we use an iterative procedure that starts from $x_0\in\mc{X}$, $z_0 \in \subdiff f(x_0, S_0)$ and generates sequences of points $x_k\in\mc{X}$ and $z_k\in\R^n$ by repeating the following steps for $k=0,1,2,\ldots$:
\begin{subequations}
\label{alg:iters:def}
\begin{align}	 
	x_{k+1} &= \argmin_{x\in\mc{X}} \left\{ \InP{z_k}{x-x_k} + \frac{1}{2\stepsize} \ltwo{x-x_k}^2\right\} \label{alg:iters:def:a}\\
	z_{k+1} &=   \mmt g_{k+1} + (1-\mmt) \frac{x_k-x_{k+1}}{\stepsize}, \label{alg:iters:def:b}	 
\end{align}
\end{subequations}
where $g_{k+1} \in \subdiff f(x_{k+1}, S_{k+1})$. When $\mc{X}=\R^n$, this algorithm reduces to the procedure \eqref{eq:SHB:unconstrained}.  
For a general convex set $\mc{X}$, this scheme is known as the iPiano method in the smooth and {deterministic} setting~\cite{OCBP14}. For simplicity, we refer to Algorithm~\ref{alg:iters:def} as stochastic heavy ball (SHB). 

Throughout the paper, we will frequently use the following two quantities:
\begin{align*}
	p_k=\frac{1-\mmt}{\mmt} \left(x_{k} - x_{k-1}\right)
	\,\,\mbox{and}\,\,\,
	d_k=\frac{1}{\stepsize} \left(x_{k-1}-x_{k}\right).
\end{align*}

Before detailing our convergence analysis, we note that most proofs of  $O(1/\epsilon^2)$ sample complexity  for subgradient-based methods rely on establishing an iterate relationship on the form (see, e.g., \cite{NJLS09,DD19,GL13}):
\begin{align}\label{eq:stoc:analysis}
	\E[V_{k+1}] \leq \E[V_k] - \alpha\, \E[e_k] + \alpha^2C^2, 
\end{align}
where $e_k$ denotes some stationarity measure such as $f(\cdot)-f\opt$ for convex and $\ltwo{\grad{f}(\cdot)}^2$ for smooth (possibly non-convex) problems, $V_k$ are certain Lyapunov functions, $\alpha$ is the stepsize, and $C$ is some constant. Once  \eqref{eq:stoc:analysis} is given, a simple manipulation results in the desired complexity provided that $\alpha$ is chosen appropriately. The case with {decaying} stepsize can be analyzed in the same way with minor adjustments. 
We follow the same route and identify a Lyapunov function that allows to establish relation \eqref{eq:stoc:analysis} for the quantity $\ltwo{\grad{F_{\lambda}(\cdot)}}$ in \eqref{eq:stationary:measure}.

\noindent Since our Lyapunov function  is nontrivial, we shall build it up through a series of key results. We start by presenting the following lemma, which quantifies the averaged progress made by one step of the algorithm. 
\begin{lemma}\label{lem:xdiff:non-smooth}
Let Assumptions (A1)--(A2) hold. Let $\mmt= \nu\stepsize$ for some constant $\nu>0$ such that $\mmt\in (0,1]$. Let $x_k$ be generated by procedure \eqref{alg:iters:def}. It holds for  any $k\in\N$ that
\begin{align}\label{eq:lem:xdiff:non-smooth}
	(1-\mmt)f(x_{k})
	+	
	\E\left[		
		\frac{\nu}{2}
		\ltwo{p_{k+1}}^2|\mc{F}_{k-1}	\right]		
	&\leq
	(1-\mmt)f(x_{k-1})
	+
	\frac{\nu}{2}
	\ltwo{p_k}^2
		- \stepsize\,\E\left[\ltwo{d_{k+1}}^2|\mc{F}_{k-1}	\right]	
	\nonumber\\
	&\hspace{0.45cm}
		+
		\stepsize^2
		\left(\frac{\rho(1-\mmt)}{2}+\nu\right)L^2.
\end{align}

\begin{proof}
See Appendix~\ref{appendix:proof:lem:xdiff:non-smooth}.
\end{proof}

\end{lemma}
In view of \eqref{eq:stoc:analysis}, the lemma shows that the quantity $\E[\ltwo{d_k}^2]$ can be made arbitrarily small. However, this alone is not sufficient to show convergence to stationary points. Nonetheless, we shall show that a small $\E[\ltwo{d_k}^2]$ indeed implies a small (averaged) value of the norm of the Moreau envelope defined at a specific point. Toward this goal, we first need to detail the points $x$ and $\hat{x}$ in \eqref{eq:stationary:measure}. It seems that taking the most natural candidate $x=x_k$ is unlikely to produce the desired result. Instead, we rely on the following iterates:
\begin{align*}
\bar{x}_{k}:= x_{k} + \frac{1-\mmt}{\mmt} \left(x_{k} - x_{k-1}\right),
\end{align*}
and construct corresponding \emph{virtual} reference points:
\begin{align*}
	\hat{x}_{k} 
	=
	\argmin_{x\in\R^n}
	\left\{
		F(x)
		+ \frac{1}{2\lambda} \ltwo{x-\bar{x}_{k}}^2\right\},
\end{align*}
for $\lambda<1/\rho$. By Lemma~\ref{lem:moreau:env}, 
$\grad{F_\lambda(\bar{x}_k)}=\lambda^{-1}(\bar{x}_k-\hat{x}_k)$,
where $F_\lambda(\cdot)$ is the Moreau envelope of $F(\cdot)=f(\cdot)+ \indcfunc{\cdot}$.

With these definitions, we can now state the next lemma.
\begin{lemma}\label{lem:Vfunc:non-smooth}
Assume the same setting of Lemma~\ref{lem:xdiff:non-smooth}. Let $\lambda>0$ be such that $\lambda^{-1}\geq 2\rho$. Let  $\xi=(1-\mmt)/\nu$ and define the function:
\begin{align}\label{eq:Vfunc}
	V_{k} 
	&= 
		F_{\lambda}(\bar{x}_{k}) 
		+
		\frac{\nu\xi^2}{4\lambda^2}
	 	\ltwo{p_{k}}^2
	 	+
	 	\frac{\stepsize\xi^2}{2\lambda^2}\ltwo{d_k}^2
	+
		\left(\frac{(1-\mmt)\xi^2}{2\lambda^2}+\frac{\xi}{\lambda}\right)  f(x_{k-1}).
\end{align}
Then, for any $k\in \N$,
\begin{align}\label{eq:lem:Vfunc:non-smooth}
	 \E\left[V_{k+1}|\mc{F}_{k-1}\right]		 
	\leq
		 V_{k}	
		-
		\frac{\stepsize}{2}		
		\ltwo{\grad{F_{\lambda}(\bar{x}_k)}}^2
		+
		\frac{\gamma\stepsize^2 L^2}{2\lambda},
\end{align}
where $\gamma = {\xi^2({\rho(1-\mmt)}/{2}+\nu)}/{\lambda} +{\rho\xi}/{2}+1$. 

\end{lemma}
The proof of this lemma is rather involved and can be found in Appendix~\ref{appendix:proof:lem:Vfunc:non-smooth}. 
Lemma~\ref{lem:Vfunc:non-smooth} has established a relation akin to \eqref{eq:stoc:analysis} with the Lyapunov function $V_k$ defined in \eqref{eq:Vfunc}. We can now use standard analysis to obtain our sample complexity.
\begin{theorem}\label{thrm:non-smooth:complexity}
Let Assumptions (A1)-(A2) hold. Let $k^*$ be sampled uniformly at random from $\{0,\ldots,K\}$. Let $f\opt=\inf_{x\in\mc{X}}f(x)$ and denote $\Delta=f(x_0)-f\opt$. If we set $\stepsize=\frac{\stepsize_0}{\sqrt{K+1}}$ and $\nu=1/\stepsize_0$ for some real $\stepsize_0>0$, then under the same setting of Lemma~\ref{lem:Vfunc:non-smooth}:
\begin{align}\label{eq:thrm:1}
	\E\left[\ltwo{\grad{F_{\lambda}(\bar{x}_{k^*})}}^2\right]	
	\leq
		2\cdot
		\frac{		
		\gamma_1\Delta
			+ 
			\frac{\gamma L^2}{2\lambda} 			
		}{
			\stepsize_0\sqrt{K+1}
		},	
\end{align}
where  $\gamma \leq \rho^2\stepsize_0^2+3\rho\stepsize_0+1$ and $\gamma_1\leq2\rho^2\stepsize_0^2+2\rho\stepsize_0+1$.
Furthermore, if $\stepsize_0$ is set to $1/\rho$, we obtain
\begin{align*}
	\E\left[\ltwo{\grad{F_{1/(2\rho)}(\bar{x}_{k^*})}}^2\right]
	\leq 
		10\cdot
	\frac{
		\rho\Delta
		+ 
		L^2
	}{
		\sqrt{K+1}
	}.	
\end{align*}
\begin{proof}
Taking the expectation on both sides of \eqref{eq:lem:Vfunc:non-smooth} and summing the result over $k=0,\ldots,K$ yield
\begin{align*}
	 \E\left[V_{K+1}\right]		 
	\leq
		 {V_{0}}
		-
		\frac{\stepsize_0}{2\sqrt{K+1}}		
		\sum_{k=0}^{K}
		\E\left[\ltwo{\grad{F_{\lambda}(\bar{x}_k)}}^2\right]
		+
		\frac{\gamma L^2\stepsize_0^2}{2\lambda}.
\end{align*}
Let $\gamma_1={1+(1-\mmt)\xi^2}/({2\lambda^2})+{\xi}/{\lambda}$, the left-hand-side of the above inequality can be lower-bounded by $\gamma_1 f\opt$. Using the facts that $F_{\lambda}(x_0)\leq f(x_0)$ and $x_{-1}=x_0$, we get $V_0\leq \gamma_1 f(x_0)$. Consequently, 
\begin{align*}
	\E\left[\ltwo{\grad{F_{\lambda}(\bar{x}_{k^*})}}^2\right]
	\leq 
		2\cdot
	\frac{
		\gamma_1\Delta
		+ 
		\frac{\gamma L^2 \stepsize_0^2}{2\lambda} 			
	}{
		\stepsize_0\sqrt{K+1}
	},	
\end{align*}
where the last expectation is taken with respect to all random sequences generated by 
the method and the uniformly distributed random variable $k^*$. 
Note that $\nu=1/\stepsize_0$, $\xi=(1-\mmt)/\nu$, and $1-\mmt\leq 1$. Thus, letting $\lambda=1/(2\rho)$, the constants  $\gamma$ and $\gamma_1$ can be upper-bounded by $\rho^2\stepsize_0^2+3\rho\stepsize_0+1$ and $2\rho^2\stepsize_0^2+2\rho\stepsize_0+1$, respectively. Therefore, if we let $\stepsize_0=1/\rho$, we arrive at
\begin{align*}
	\E\left[\ltwo{\grad{F_{1/(2\rho)}(\bar{x}_{k^*})}}^2\right]
	\leq 
		10\cdot
	\frac{
		\rho\Delta
		+ 
		L^2
	}{
		\sqrt{K+1}
	},	
\end{align*}
as desired.
\end{proof}

\end{theorem}

Some remarks regarding Theorem~\ref{thrm:non-smooth:complexity} are in order: 

i) The choice $\nu=1/\stepsize_0$ is just for simplicity; one can pick any constant such that $\mmt=\nu\stepsize\in(0,1]$. 
Note that the stepsize used to achieve the rate in \eqref{eq:thrm:1} does not depend on any problem parameters. Once $\stepsize$ is set, the momentum parameter selection is completely parameter-free. Since both $\stepsize$ and $\beta$ scale like $O(1/\sqrt{K})$, Algorithm~\ref{alg:iters:def} can be seen as a single time-scale method \cite{GRW18, RS83}. Such methods contrast those that require at least two time-scales to ensure convergence. For example, stochastic dual averaging for \emph{convex} optimization \cite{Xia10} requires one fast scale $O(1/K)$ for averaging the subgradients, and one slower scale $O(1/\sqrt{K})$ for the stepsize. To show almost sure convergence of SHB for \emph{smooth and unconstrained} problems, the work \cite{GLZX19} requires that both the stepsize and the momentum parameter tend to zero but the former one must do so at a faster speed.\footnote{Note that our $\mmt$ corresponds to $1-\mmt$ in \cite{GLZX19}.}  \medskip

ii) To some extent, Theorem~\ref{thrm:non-smooth:complexity} supports the use of a small momentum parameter such as $\mmt=0.1$ or $\mmt=0.01$, which corresponds to the default value $1-\mmt=0.9$ in PyTorch\footnote{\url{https://pytorch.org/}} or the smaller $1-\mmt=0.99$ suggested in~\cite{Goh17}.
Indeed, the theorem allows to have $\mmt$ as small as $O(1/\sqrt{K})$, i.e., one can put much more weight to the momentum term than the fresh subgradient and still preserve the complexity. Recall also that SHB reduces to SGD as $\mmt\to 1$, which also admits a similar complexity. It is thus quite flexible to set $\mmt$, without sacrificing the worst-case complexity. We refer to \cite[Theorem~2]{GLZX19} for a similar discussion in the context of almost sure convergence on smooth problems. 

iii) In view of \eqref{eq:stationary:measure}, the theorem indicates that $\bar{x}_k$ is nearby a near-stationary point $\hat{x}_k$. Since $\bar{x}_k$ may not belong to $\mc{X}$, it is thus more preferable to have the similar guarantee for the iterate $x_k$. Indeed, we have
\begin{align*}
	\lambda^{-2}\ltwo{x_k-\hat{x}_k} ^2
	& \leq  
		 2\lambda^{-2}\ltwo{\bar{x}_k-\hat{x}_k}^2  + 2\lambda^{-2}\ltwo{x_k-\bar{x}_k}^2
	\\
	&=
		2\ltwo{\grad{F_{\lambda}(\bar{x}_{k})}}^2+2\lambda^{-2}\ltwo{x_k-\bar{x}_k}^2
	\\
	&=
		2\ltwo{\grad{F_{\lambda}(\bar{x}_{k})}}^2+2\lambda^{-2}\xi^2\ltwo{d_k}^2.
\end{align*}
Since both terms on the right converge at the rate $O(1/\sqrt{K})$, it immediately translates into the same guarantee for the term on the left, as desired. 

In summary, we have established the convergence rate $O(1/\sqrt{K})$ or, equivalently, the sample complexity $O(1/\epsilon^2)$ of SHB  for the minimization of weakly convex functions. 

\section{Extension to smooth non-convex functions}

In this section, we study the convergence property of Algorithm~\eqref{alg:iters:def} for the minimization of $\rho$-smooth functions:
\begin{align*}
	\ltwo{\grad{f(x)}-\grad{f(x)}}\leq\rho\ltwo{x-y}, \,\,\, \forall  x,y\in \dom f.
\end{align*}
Note that $\rho$-smooth functions are automatically $\rho$-weakly convex.
In this setting, it is more common to replace Assumption (A2) by the following.

\noindent\textbf{Assumption (A3)}.  There exists a real $\sigma>0$ such that  for all $x\in \mc{X}$:
\begin{align*}
	\E\left[\ltwo{f'(x, \statrv) - \grad{f}(x)}^2\right] \leq \sigma^2.
\end{align*}
Deriving convergence rates of stochastic schemes with momentum for non-convex functions under Assumption~(A3) can be quite challenging. Indeed, even in \emph{unconstrained} optimization, previous studies often need to make the assumption that the true gradient is bounded, i.e., $\ltwo{\grad{f(x)}}\leq G$ for all $x\in\R^n$ (see, e.g., \cite{YYL18,GLZX19}). This assumption is strong and does not hold even for \emph{quadratic convex} functions. It is more realistic in constrained problems,  for example when $\mc{X}$ is compact, albeit the constant $G$ could then be  large. 

Our objective in this section is twofold: First, we aim to extend the convergence results of SHB in the previous section to smooth optimization problems under Assumption~(A3). Note that even in this setting, the sample complexity of SHB has not been established before. The rate is obtained without the need of forming a batch of samples at each iteration, which is commonly required for constrained non-convex stochastic optimization \cite{GL16, GLZ16}. Second, for unconstrained problems, we demonstrate how to achieve the same complexity without the bounded gradient assumption above.

Let $h(x)=\frac{1}{2}\ltwo{x}^2+\indcfunc{x}$ and let $h^*(z)$ be its convex conjugate. Our convergence analysis relies on the  function:
\begin{align}
	\varphi_k 
	= 
		h^*(x_{k} - \stepsize z_{k}) 
		- 
		\half \ltwo{x_{k}}^2 + \stepsize \InP{x_{k}}{ z_{k}}.
\end{align}
The use of this function is inspired by \cite{Rus87}.  Roughly speaking, $\varphi_k$ is the negative of the optimal value of the function on the RHS of \eqref{alg:iters:def:a}, and hence, $\varphi_{k}\geq 0$ for all $k$.  
This function also underpins the analysis of  the dual averaging scheme in \cite{Nes09}.

The following result plays a similar role as Lemma~\ref{lem:xdiff:non-smooth}. 
\begin{lemma}\label{lem:Wfunc:smooth}
Let  Assumptions~(A1) and (A3) hold. Let $\stepsize\in(0, 1/\rho)$ and $\mmt= \nu\stepsize$ for some constant $\nu>0$ such that $\mmt\in (0,1]$. 
Let $\stepsize\in\left(0,1/(4\rho)\right]$ and $\xi=(1-\mmt)/\nu$, and define the function:
\begin{align*}
	W_k = 2f(x_k) + \frac{\varphi_k}{\nu\stepsize^2} + \frac{\xi}{2}\ltwo{d_k}^2.
\end{align*}
Then, it holds for any $k\in\N$ that
\begin{align}\label{eq:lem:smooth:W}
	\E\left [W_{k+1} | \mc{F}_{k}\right ]
	&\leq
		W_k
		-\stepsize
		\ltwo{d_{k+1}}^2
		+
		4\nu\stepsize^2\sigma^2.
\end{align}
\begin{proof}
Since the proof is rather technical,  we defer details to Appendix~\ref{appendix:proof:Wfunc:smooth} and sketch only the main arguments here. 

By smoothness  of $h^*$, weak convexity of $f$ and the optimality condition for the update formula \eqref{alg:iters:def:a}  we get
\begin{align*}
	\E\left[f(x_{k+1}) + \frac{\varphi_{k+1}}{\nu\stepsize^2} \Big| \mc{F}_{k}\right ]
	&\leq
		 f(x_k)+ \frac{\varphi_{k}}{\nu\stepsize^2}
		-(\stepsize -\frac{\rho\stepsize^2 }{2})
		\ltwo{d_{k+1}}^2
	 	+
	 	\frac{1}{2\nu}\E\left [ \ltwo{ z_{k} -  z_{k+1}}^2| \mc{F}_{k}\right ].
\end{align*}
The proof of this relation can be found in Lemma~\ref{lem:xdiff:smooth}.
The preceding inequality admits very useful properties as we have terms that form a telescoping sum, and the constant associated with $\ltwo{d_{k+1}}^2$ has the right order-dependence on the stepsize. However, we still have a remaining term $\ltwo{ z_{k} -  z_{k+1}}^2$. Thus, in view of relation \eqref{eq:stoc:analysis}, our next strategy  is to bound this term in a way that still keeps all the favourable features described above, and at the most introduces an additional term of order $O(\stepsize^2 \sigma^2)$. As shown in Lemma~\ref{lem:z:diff}, we can establish the following inequality
\begin{align*}
	\E\left[\frac{1}{2\nu}\ltwo{z_{k+1}-z_k}^2|\mc{F}_k\right]
	&\leq
		f(x_k)-f(x_{k+1})
		+		
		\frac{\xi}{2}\ltwo{d_k}^2	- \frac{\xi}{2} \ltwo{d_{k+1}}^2	
	\nonumber\\	 	
	&\hspace{0.45cm}
		-\left(\stepsize-\frac{\stepsize^3\rho^2+3\rho\stepsize^2}{2}\right)\ltwo{d_{k+1}}^2
		+
		4\nu\stepsize^2\sigma^2.
\end{align*} 
Now, (\ref{eq:lem:smooth:W}) follows immediately from combining the two previous inequalities and  the fact that $\stepsize\in\left(0,1/(4\rho)\right]$. 
\end{proof}
\end{lemma}

We remark that Lemma~\ref{lem:Wfunc:smooth} does not require the bounded gradient assumption and readily indicates the convergence rate $O(1/\sqrt{K})$ for $\E[\ltwo{d_k}^2]$. However,  to establish the rate for $\E[\ltwo{\grad{F_{\lambda}(\bar{x}_k)}}^2]$, we  need to impose such an assumption in the theorem below. Nonetheless, the assumption is much more realistic in this setting than the unconstrained case.  
\begin{theorem}\label{thrm:smooth:complexity}
Let Assumptions (A1) and (A3) hold. Assume further that $\ltwo{\grad{f(x)}}\leq G$ for all $x\in \mc{X}$. Let $k^*$, $\bar{x}_{k^*}$, $\lambda$, $\Delta$, $\gamma$, and $\gamma_1$ be defined as in Theorem~\ref{thrm:non-smooth:complexity}. If we set $\stepsize=\frac{\stepsize_0}{\sqrt{K+1}}$ and $\nu=1/\stepsize_0$ for some real $\stepsize_0>0$, then 
\begin{align*}
	\E\left[\ltwo{\grad{F_{\lambda}(\bar{x}_{k^*})}}^2\right]	
	\leq
		2\cdot
		\frac{		
		\gamma_1\Delta
			+ 
			{\gamma (\sigma^2+G^2)}/({2\lambda} )			
		}{
			\stepsize_0\sqrt{K+1}
		}.	
\end{align*}
Furthermore, if $\stepsize_0$ is set to $1/\rho$, we obtain
\begin{align*}
	\E\left[\ltwo{\grad{F_{1/(2\rho)}(\bar{x}_{k^*})}}^2\right]
	\leq 
		10\cdot
	\frac{
		\rho\Delta
		+ 
		\sigma^2+G^2
	}{
		\sqrt{K+1}
	}.	
\end{align*}
\begin{proof}
The proof is a verbatim copy of that of Theorem~\ref{thrm:non-smooth:complexity} with $L^2$ replaced by $\sigma^2+G^2$; see Appendix~\ref{appendix:proof:thrm:smooth:complexity}. 
\end{proof}
\end{theorem}
Some remarks are in order: 

i) To the best of our knowledge, this is the first convergence rate result of a stochastic (or even deterministic) method with Polyak momentum for smooth, non-convex, and constrained problems. 

ii) The algorithm enjoys the same single time-scale and parameter-free properties as in the non-smooth case.

iii) The rate in the theorem readily translates into an analogous estimate for  the norm of the so-called \emph{gradient mapping} $\mc{G}_{1/\rho}$, which is commonly adapted in the literature, e.g., \cite{GLZ16}. This is because for $\rho$-smooth functions \cite{DP19}:
\begin{align*}
		\ltwo{\mc{G}_{1/\rho}(x)}
	\leq 
		\frac{3}{2}(1+ \frac{1}{\sqrt{2}})\ltwo{\grad{F_{1/(2\rho)}(x)}}, \,\,\, \forall x \in \R^n.
\end{align*}


Since the bounded gradient assumption is rather restrictive in the  unconstrained  case, our final result shows how the desired complexity can be attained without this assumption. 
\begin{theorem}\label{thrm:smooth:unconstrained:complexity}
Let Assumptions (A1) and (A3) hold. Let $\lambda^{-1}\in(3\rho/2,2\rho]$.  Let $k^*$ be sampled uniformly at random from $\{-1,\ldots,K-1\}$. If we set $\stepsize=\frac{\stepsize_0}{\sqrt{K+1}}$ and $\nu=1/\stepsize_0$, where $\stepsize_0\in(0,1/(4\rho)]$, then under the same setting of Lemma~\ref{lem:Wfunc:smooth}:
\begin{align*}
	\E\left[\ltwo{\grad{F_{\lambda}(\bar{x}_{k^*})}}^2\right]
	\leq
		c\cdot
	\frac{
		(1 +2\stepsize_0^2/\lambda^2)\Delta
		+ 
		\frac{(1+{8\stepsize_0}/{\lambda}) \sigma^2\stepsize_0^2}{2\lambda}	
	}{
		\stepsize_0\sqrt{K+1}
	},
\end{align*}
where $c={2\lambda^{-1}}/{(2\lambda^{-1}-3\rho)}$.
Furthermore, let $\lambda=1/(2\rho)$, we obtain
\begin{align*}
	&\E\left[\ltwo{\grad{F_{1/(2\rho)}(\bar{x}_{k^*})}}^2\right]	
	\leq
		4\cdot
	 	\frac{
	 		\left(1 +{8\rho^2\stepsize_0^2}\right)\Delta
	 		+ 
	 		{(\rho+16\stepsize_0\rho^2)}\sigma^2\stepsize_0^3	
	 	}{
	 		\stepsize_0\sqrt{K+1}
	 	}.
\end{align*}
\begin{proof}
See Appendix~\ref{appendix:proof:thrm:smooth:unconstrained:complexity}.
\end{proof}
\end{theorem}

It should be mentioned that a similar result has been attained very recently in \cite{GRW18} using a different analysis, albeit no sample complexity is given for the non-smooth case. It is still an open question if one can preserve the complexity in Theorem~\ref{thrm:smooth:complexity} without the bounded gradient hypothesis.

\section{Numerical evaluations}
In this section, we perform experiments to validate our theoretical developments and to demonstrate that despite sharing the same worst-case complexity, SHB can be better in terms of speed and robustness to problem and algorithm parameters than SGD.

We consider the robust phase retrieval problem \cite{DR18a,DR18b}: Given a set of $m$ measurements $(a_i, b_i) \in \R^{n}\times \R$, the phase retrieval problem seeks for a vector $x\opt$ such that $\InP{a_i}{x\opt}^2 \approx b_i$ for most $i=1,\ldots,m$.
Whenever the problem is corrupted with gross outliers, a natural exact penalty form of this (approximate) system of equations yields the minimization problem: 
\begin{align*}
	\minimize_{x\in\R^n} \frac{1}{m} \sum_{i=1}^{m} \big|\InP{a_i}{x}^2-b_i \big|.
\end{align*}
This objective function is non-smooth and non-convex. In view of \eqref{eq:composite}, it is the composition of the Lipschitz-continuous convex function $h(y)=\norm{y}_1$ and the smooth map $c$ with $c_i(x)=\InP{a_i}{x}^2-b_i$\replaced{. Hence, it is weakly convex.}{,thereby being weakly convex.}
\begin{figure}[t!]
	\centering
	\begin{minipage}{0.495\textwidth}
		\centering
		{\includegraphics[width=1.\textwidth]{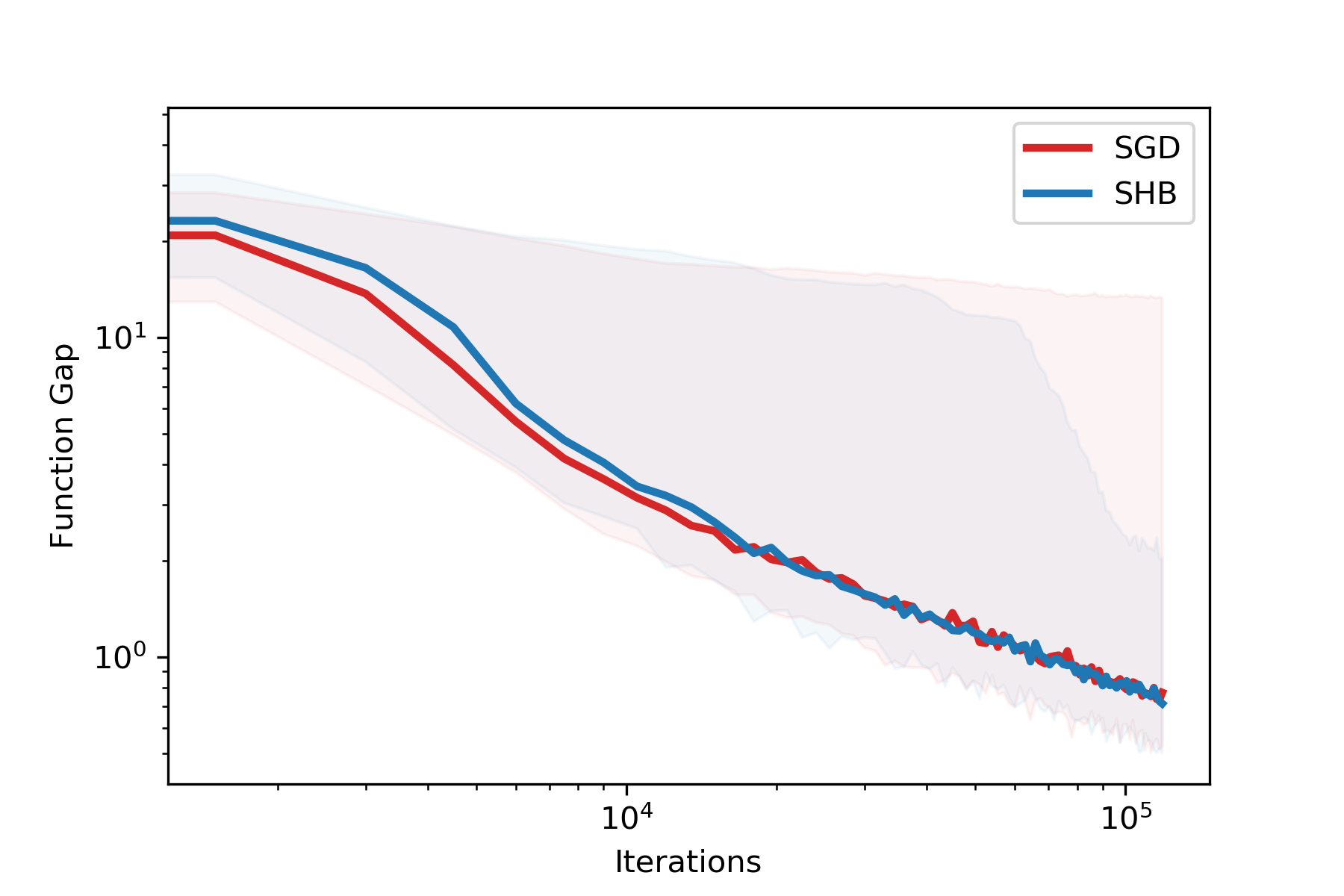}}
		\subcaption{$\kappa=10$, $\stepsize_0=0.1$}	
	\end{minipage}
\hskip -0.15in
	\begin{minipage}{0.495\textwidth}
		\centering
		{\includegraphics[width=1.\textwidth]{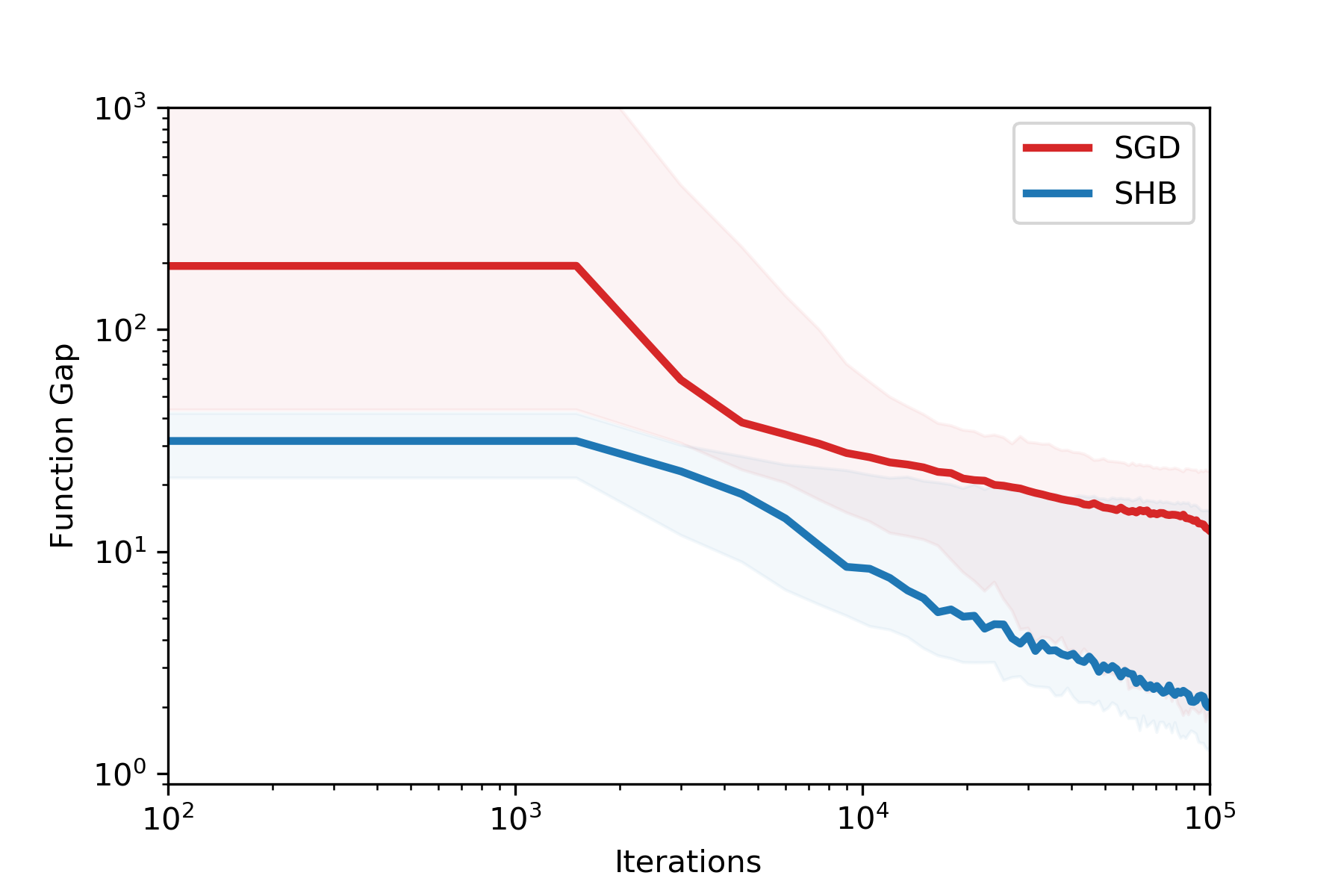}}
		\subcaption{$\kappa=10$, $\stepsize_0=0.25$}	
		\end{minipage}
	\begin{minipage}{0.495\textwidth}
		\centering 
		{\includegraphics[width=1.\textwidth]{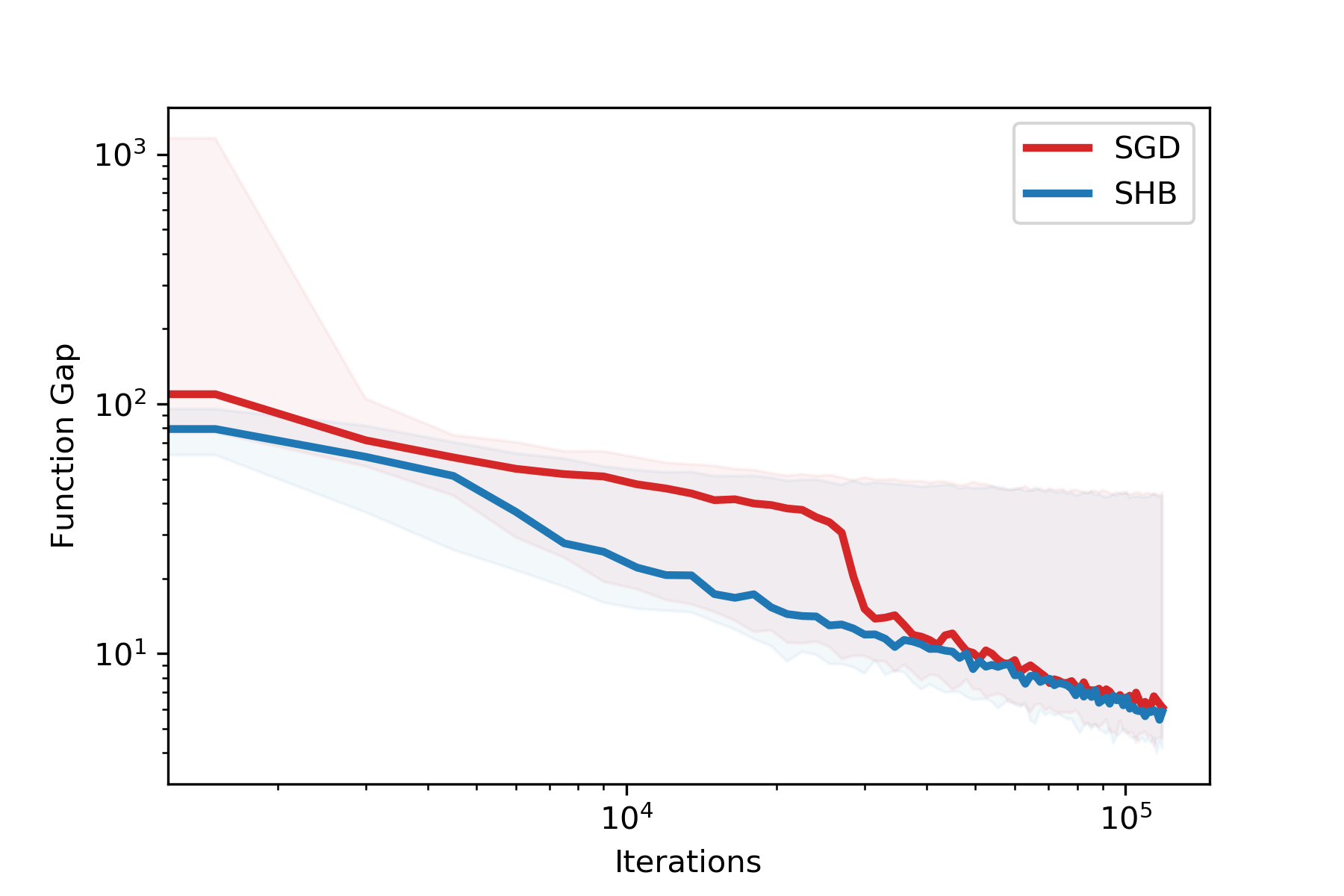}}
		\subcaption{$\kappa=1$, $\stepsize_0=0.1$}	
	\end{minipage}
\hskip -0.15in
	\begin{minipage}{0.495\textwidth}
	\centering
	{\includegraphics[width=1.\textwidth]{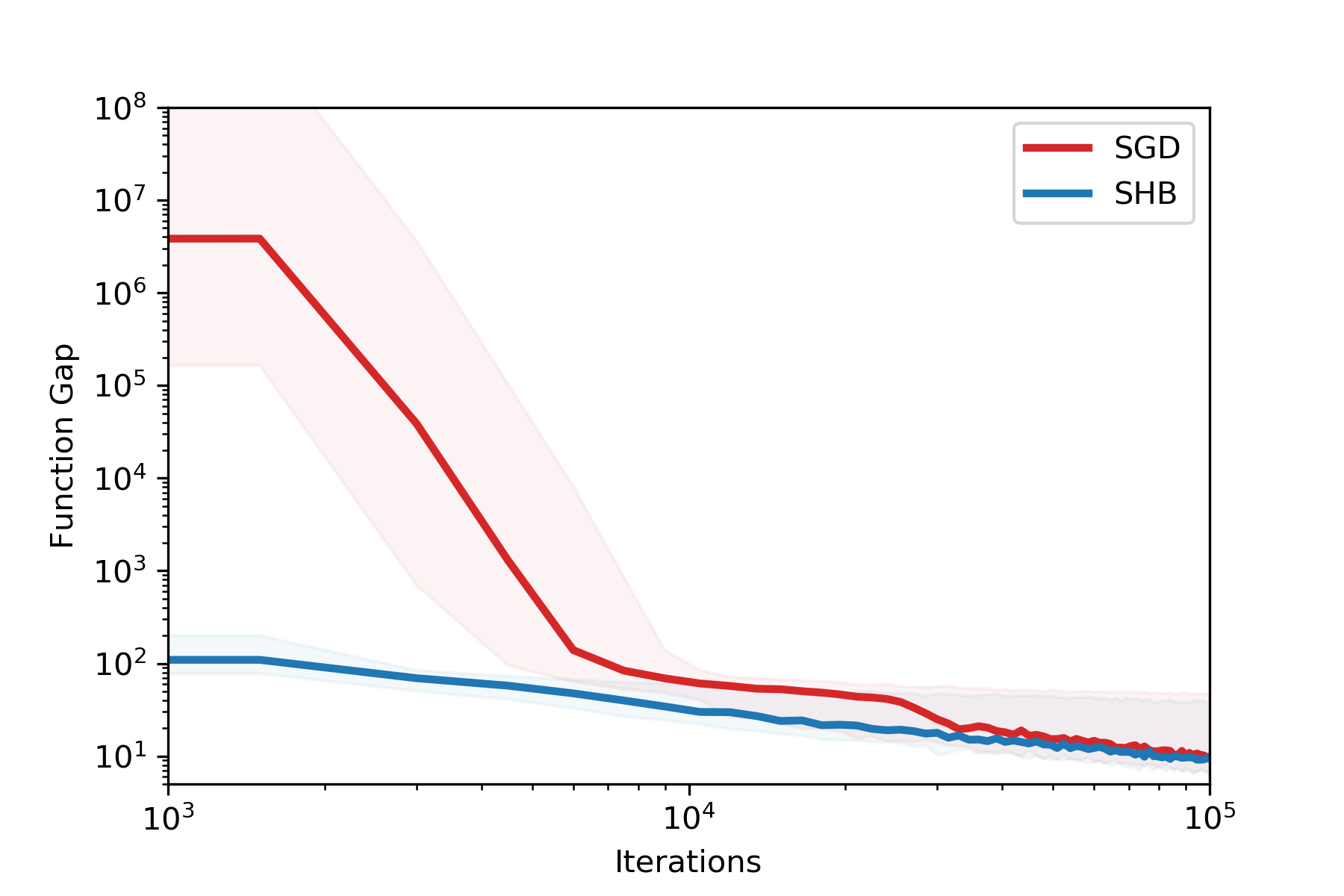}}
	\subcaption{$\kappa=1$, $\stepsize_0=0.15$}	\label{fig:1:d}
	\end{minipage}
	\caption{The function gap $f(x_k)-f(x\opt)$ versus iteration count for phase retrieval with $p_\textrm{fail}=0.2$, $\mmt=10/\sqrt{K}$. For better visualization,  in this plot we generated  $x^\star$ with standard normal distributed elements.}\label{fig:1}
\end{figure}

\begin{figure*}[!t]
	\centering
	\begin{minipage}{0.495\textwidth}
		\centering
		{\includegraphics[width=1.\textwidth]{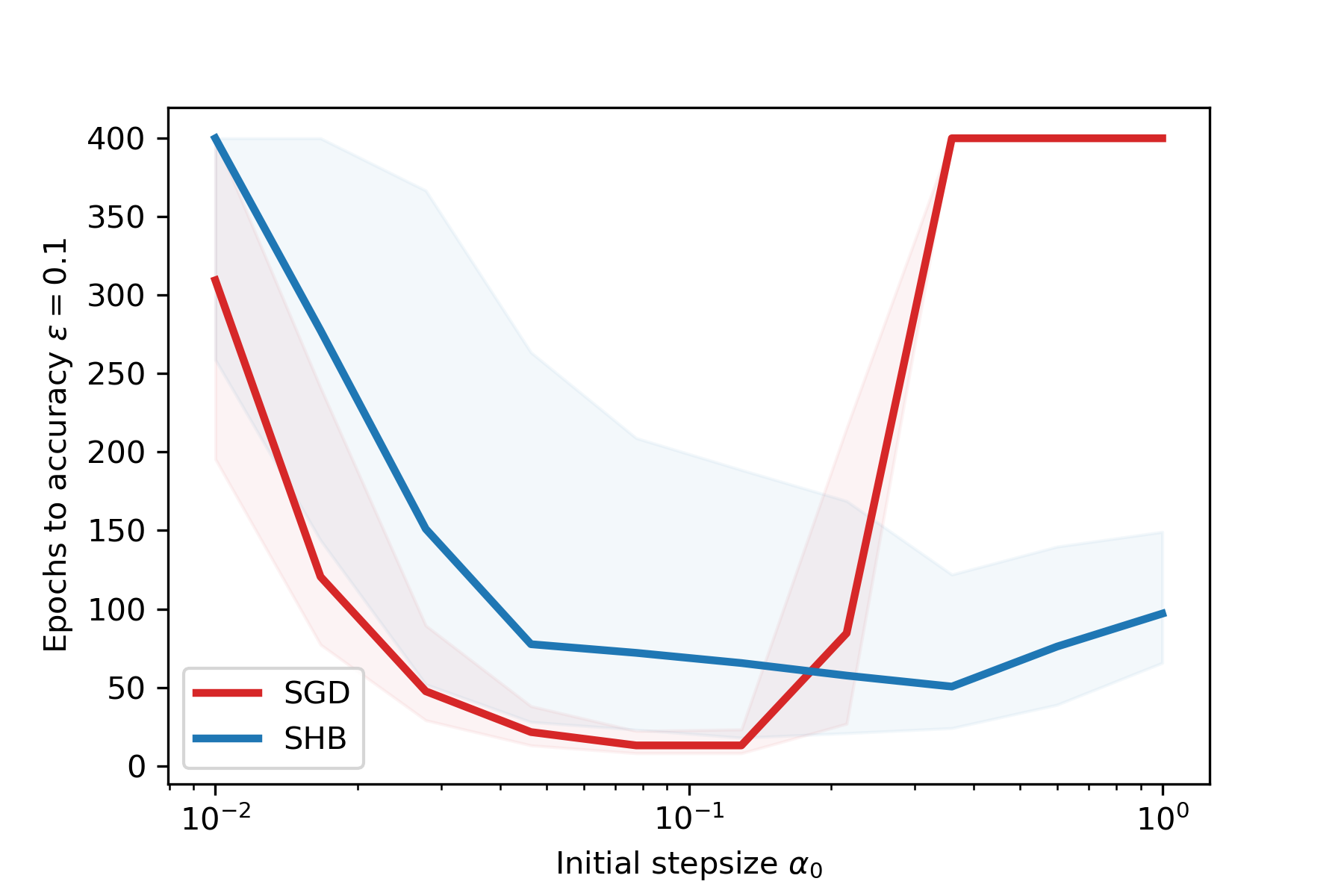}}
		\label{fig:epochs2eps} 
		\subcaption{ $p_\textrm{fail}=0.3$, $\mmt=1/\sqrt{K}$}
	\end{minipage}
	~
	\begin{minipage}{0.45\textwidth}
		\centering
		{\includegraphics[width=1.\textwidth]{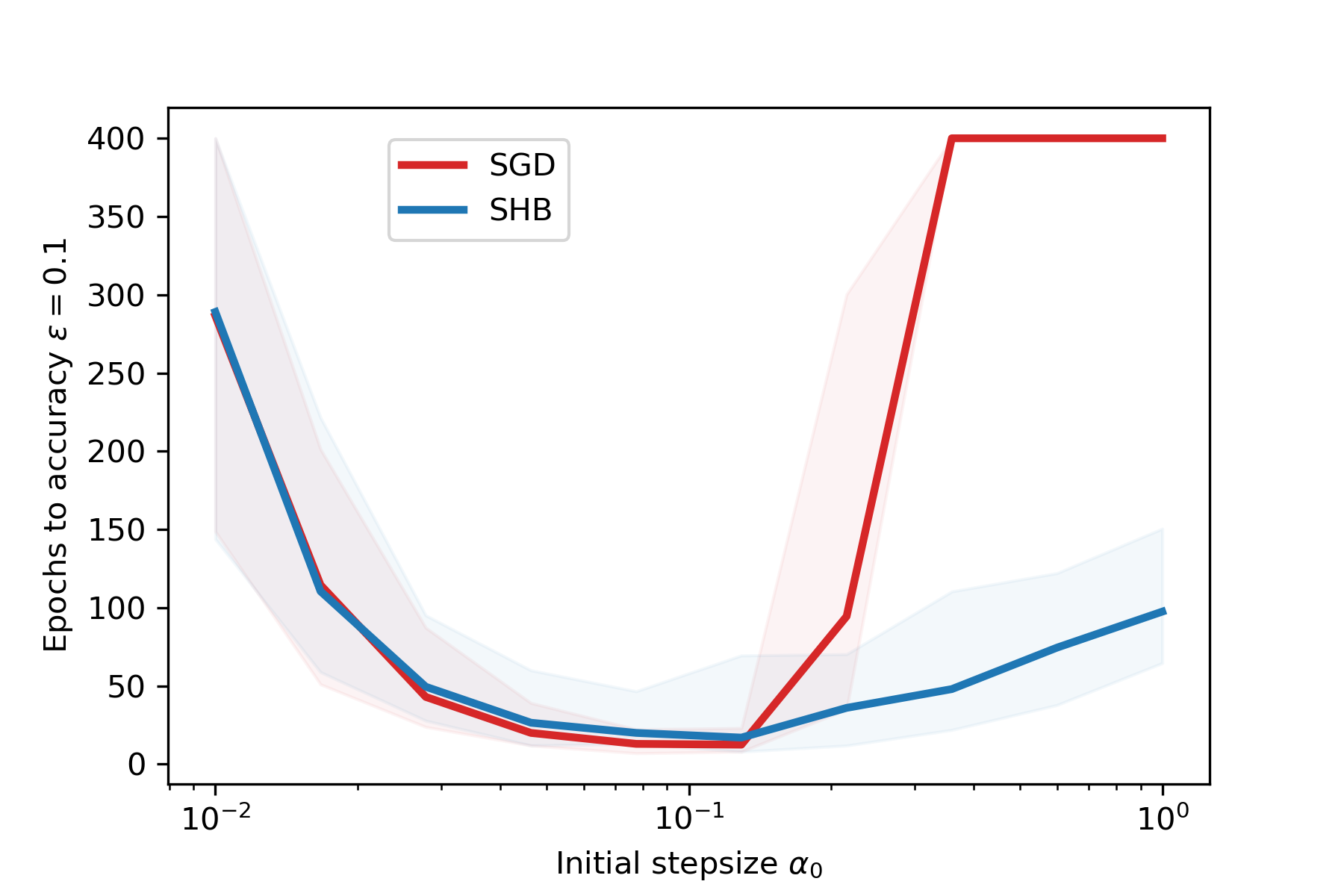}}
		\label{fig:epochs2eps:corrupted} 
		\subcaption{$p_\textrm{fail}=0.3$, $\mmt=1/\stepsize_0/\sqrt{K}$}
	\end{minipage}
	\caption{The number of epochs to achieve $\epsilon$-accuracy versus initial stepsize $\stepsize_0$ for phase retrieval with $\kappa=10$.}\label{fig:2}
\end{figure*}
\begin{figure*}[!t]
	\centering
	\begin{minipage}{0.45\textwidth}
		\centering
		{\includegraphics[width=1.\textwidth]{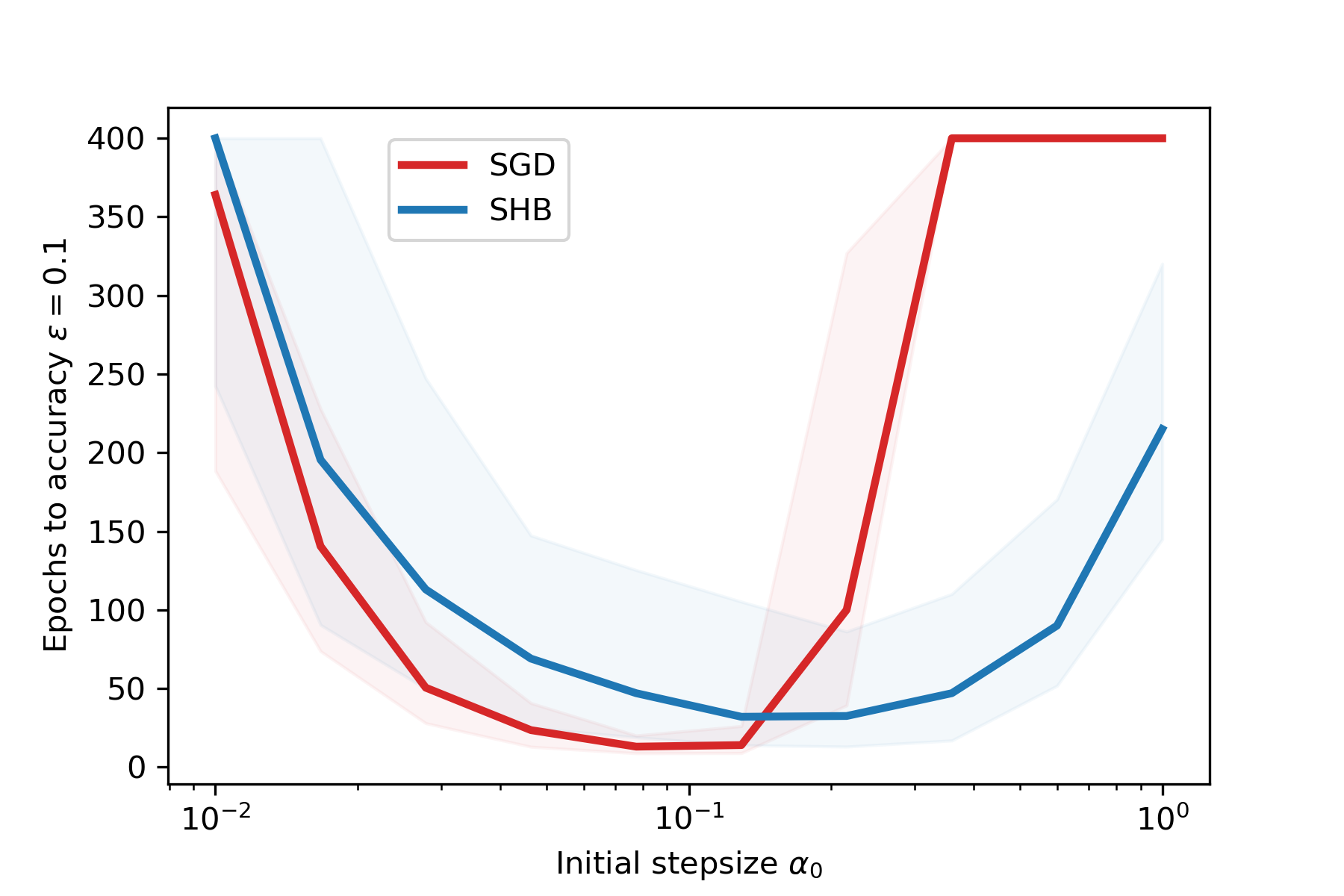}}
		\label{fig:epochs2eps:popular:mmt} 
		\subcaption{$p_\textrm{fail}=0.3$, $\mmt=0.1$}
	\end{minipage}
	~
	\begin{minipage}{0.45\textwidth}
		\centering
		{\includegraphics[width=1.\textwidth]{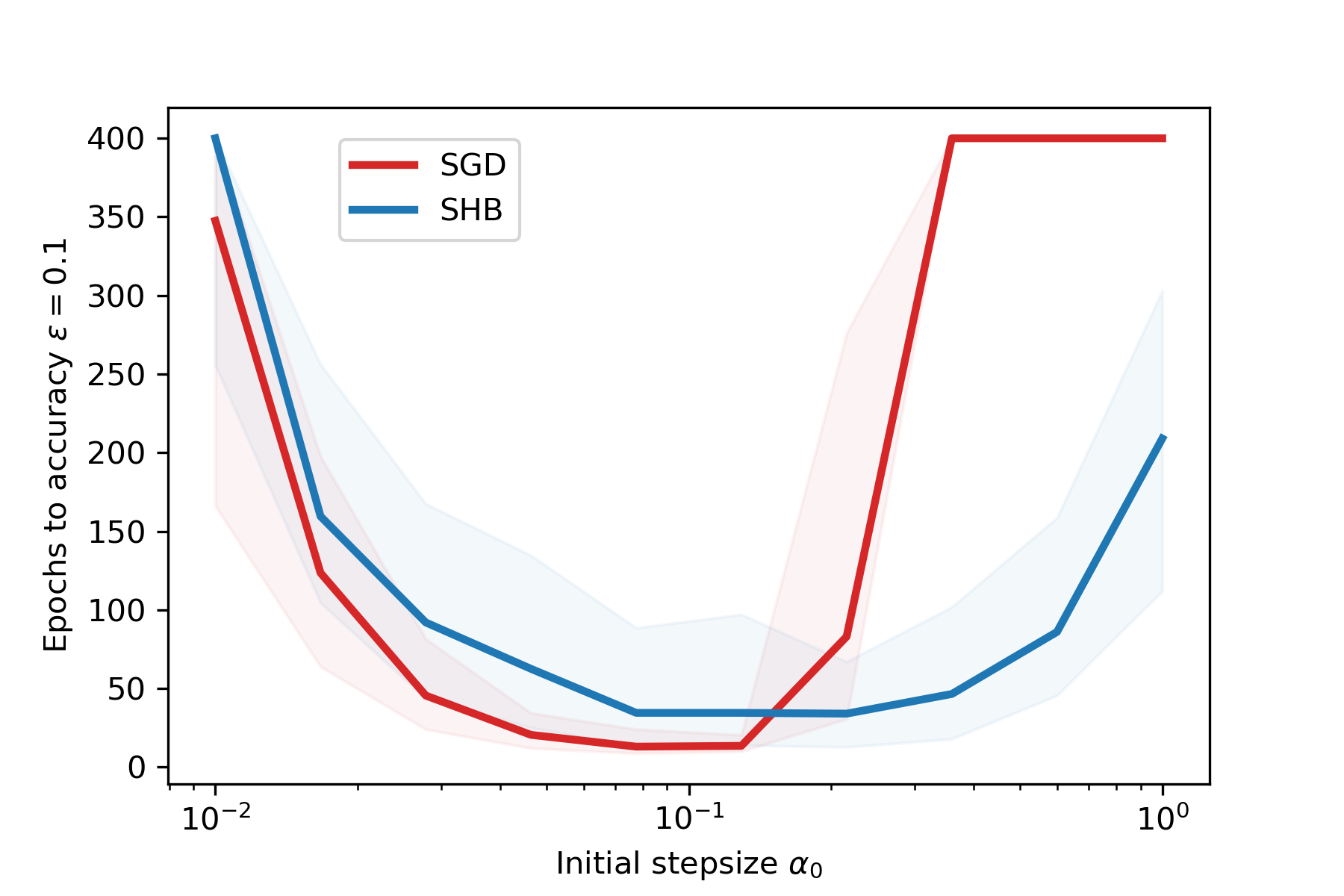}}
		\label{fig:epochs2eps:corrupted:popular:mmt} 
		\subcaption{$p_\textrm{fail}=0.3$, $\mmt=0.01$}
	\end{minipage}
	\caption{The number of epochs to achieve $\epsilon$-accuracy versus initial stepsize $\stepsize_0$ for phase retrieval with $\kappa=10$ and popular choices of $\mmt$.}\label{fig:3}
\end{figure*}

In each experiment, we set $m=300$\replaced{, }{ and} $n=100$\deleted{,}  and select $x\opt$ uniformly from the unit sphere. We  generate $A$ as $A=Q D$, where $Q\in\R^{m\times n}$ is a matrix with standard normal distributed entries, and $D$ is a diagonal matrix with linearly spaced elements between $1/\kappa$ and $1$\added{, } with $\kappa\geq 1$ playing the role of a condition number. The elements $b_i$ of the vector $b$ are generated as $b_i = \InP{a_i}{x\opt}^2 + \delta \zeta_i$,  $i =1, \ldots,m,$ where $\zeta_i \sim \mc{N}(0, 25)$ models the corruptions, and $\delta\in\{0,1\}$ is a binary random variable taking the value $1$ with probability $p_{\textrm{fail}}$, so that $p_{\textrm{fail}} \cdot m$ measurements \replaced{are corrupted}{provide no information}. The algorithms are all 
randomly initialized at $x_0\sim \mc{N}(0, 1)$. The stochastic subgradient is simply given as an element of the subdifferential of $g(x)=\big|\InP{a}{x}^2-b \big|$:
\begin{align*}
	\subdiff g(x) 
	= 
		2\InP{a}{x}a
		\cdot
		\begin{cases}
			\sign{(\InP{a}{x}-b)}^2 	& \mbox{if}  \,\,  \InP{a}{x}^2 \neq b,\\
			[-1,1] 								&   \mbox{otherwise}.
		\end{cases}
\end{align*}

In each of our experiments, we set the stepsize as $\stepsize_k=\stepsize_0/\sqrt{k+1}$, where $\stepsize_0$ is an initial stepsize. We note that this stepsize can often make a little faster progress (for both SGD and SHB) at the beginning of the optimization process than the constant one $\stepsize_k=\stepsize_0/\sqrt{K+1}$. However, after a few iterations, both of them yield very similar results and will not change the qualitative aspects of our plots in any way. We also refer $m$ stochastic iterations as one epoch (pass over \added{the} data). Within each individual run, we allow the considered stochastic methods to perform $K=400m$ iterations. 
We conduct 50 experiments for each stepsize and report the median of the so-called epochs-to-$\epsilon$-accuracy; the shaded areas in each plot cover the 10th to 90th percentile of convergence times. Here, the  epoch-to-$\epsilon$-accuracy is defined as the smallest number of epochs $q$ required to reach $f(x_{m\cdot  q}) -  f(x^\star) \leq \epsilon$.

Figure~\ref{fig:1} shows the function gap versus iteration count for different values of $\kappa$ and $\stepsize_0$, with $p_\textrm{fail}=0.2$, $\mmt=10/\sqrt{K}$. It is evident that SHB converges with a theoretically justified parameter $\mmt$ and is much less sensitive to problem and algorithm parameters than the vanilla SGD. Note that the sensitivity issue of SGD is rather well documented; a \replaced{slight}{slightly} change in \replaced{its}{such} parameters can have a severe effect on the overall performance of the algorithm \cite{NJLS09, AD19a}. For example, Fig.~\ref{fig:1:d} shows that SGD  exhibits a transient exponential growth before eventual convergence. This behaviour \replaced{can occur even when}{persists even for} minimizing the smooth quadratic function $\half x^2$ \cite[Example~2]{AD19a}. In contrast, SHB converges in all settings of the figure, suggesting that using a momentum term can help to improve the robustness of the standard SGD.  This is expected as the update formula \eqref{eq:SHB:unconstrained:z} acts like a \added{lowpass} filter, averaging past stochastic subgradients, which may have stabilizing effect on the sequence $\{x_k\}$.

To further clarify this observation, in the next set of experiments, we test the sensitivity of SHB and SGD to the initial stepsize $\stepsize_0$. Figure~\ref{fig:2} shows the number of epochs required to reach $\epsilon$-accuracy for phase retrieval with $\kappa=10$ and $p_\textrm{fail}=0.3$. 
We can see that the standard SGD has good performance for a narrow range of stepsizes, while wider convergence range can be achieved with SHB. 

Finally, it would be incomplete without reporting experiments for some of the most popular momentum parameters used by practitioners. Figure~\ref{fig:3} shows a similar story to Fig.~\ref{fig:2} for the parameter $1-\mmt=0.9$ and $1-\mmt=0.99$ as discussed in remark ii) after Theorem~\ref{thrm:non-smooth:complexity}. 
This together with Fig.~\ref{fig:2} demonstrates that SHB is able to find good approximate solutions for diverse values of the momentum constant over a wider (often significantly so) range of algorithm parameters than SGD.

\section{Conclusion}

Using a carefully constructed Lyapunov function, we established the first sample complexity results for the SHB method on a broad class of non-smooth, non-convex, and constrained optimization problems. The complexity is attained in a parameter-free fashion in a single time-scale. A notable feature of our results is that they justify the use of a large amount of momentum in search directions. We also improved some complexity results for SHB on smooth problems. Numerical results show that SHB exhibits good performance and low sensitivity to problem and algorithm parameters compared to the standard SGD.

\section*{Acknowledgements}

This work was supported in part by the Knut and Alice Wallenberg Foundation, the Swedish Research Council and the Swedish Foundation for Strategic Research. We would like to thank Ahmet Alacaoglu for his useful feedback on an earlier version of this manuscript.
\nocite{GFJ15,ZK93}

\bibliographystyle{abbrv}
\bibliography{refs}

\appendix

\section{Proof of Lemma~\ref{lem:xdiff:non-smooth}}\label{appendix:proof:lem:xdiff:non-smooth}

First, we write the update formula of $x_{k+1}$ in \eqref{alg:iters:def} on a more compact form as:
\begin{align}
	x_{k+1} &= \proj{x_k-\stepsize z_k}\\
	z_{k+1} &= \mmt g_{k+1} +(1-\mmt)\frac{x_k-x_{k+1}}{\stepsize}.
\end{align}
We have
\begin{align*}
	\frac{1}{2}
	\ltwo{\frac{1-\mmt}{\mmt} \left(x_{k+1} - x_{k}\right)}^2			
	=
		\frac{1}{2\mmt^2}
		\ltwo{x_{k+1} -  \left((1-\mmt)x_{k} + \mmt x_{k+1} \right)}^2.
\end{align*}
Since $(1-\mmt)x_{k} + \mmt x_{k+1} \in\mc{X}$ and $\proj{\cdot}$ is nonexpansive, it holds that
\begin{align}\label{eq:proof:lem:1:norm}
	\frac{1}{2}
	\ltwo{\frac{1-\mmt}{\mmt} \left(x_{k+1} - x_{k}\right)}^2			
	&\leq
		\frac{1}{2\mmt^2}
		\ltwo{x_k-\stepsize z_k -  \left((1-\mmt)x_{k} + \mmt x_{k+1}\right)}^2		
	\nonumber\\
	&=
		\frac{1}{2}
		\ltwo{ x_{k} + \frac{1-\mmt}{\mmt} \left(x_{k} - x_{k-1}\right)
		- \stepsize g_k 
		-x_{k+1}}^2.	
\end{align}
Next, we decompose the right-hand-side of the preceding inequality as 
\begin{align}\label{eq:proof:lem:1:norm:decomp}
	\frac{1}{2}
	\ltwo{ x_{k} + \frac{1-\mmt}{\mmt} \left(x_{k} - x_{k-1}\right) - \stepsize g_k -x_{k+1}}^2
	&=
		\frac{1}{2}
		\ltwo{\frac{1-\mmt}{\mmt} \left(x_{k} - x_{k-1}\right)}^2	
	\nonumber\\
	&\hspace{0.45cm}
		+
		\frac{(1-\mmt)\stepsize}{\mmt}\InP{g_k}	{x_{k-1}-x_k}
	\nonumber\\
	&\hspace{0.45cm}
		+
		\frac{1-\mmt}{\mmt}\InP{x_k-x_{k-1}}{x_k-x_{k+1}}	
		+			
		\stepsize\InP{g_k}{x_{k+1}-x_k}	
	\nonumber\\
	&\hspace{0.45cm}
		+
		\frac{1}{2}\ltwo{x_{k+1}-x_k}^2 
		+
		\frac{\stepsize ^2}{2}\ltwo{g_k}^2.
\end{align}
Let $p_k:=\frac{1-\mmt}{\mmt} \left(x_{k} - x_{k-1}\right)$, then by multiplying both sides of \eqref{eq:proof:lem:1:norm} by $\nu=\mmt/\stepsize>0$, together with \eqref{eq:proof:lem:1:norm:decomp} and  the definition of $z_k$, we get
\begin{align}\label{eq:proof:lem:1:norm:2}
	\frac{\nu}{2}
	\ltwo{p_{k+1}}^2
	&\leq
		\frac{\nu}{2}
		\ltwo{p_k}^2	
		+
		(1-\mmt)\InP{g_k}	{x_{k-1}-x_k}
		+
		\InP{z_k}{x_{k+1}-x_k}	
	\nonumber\\
	&\hspace{0.45cm}
		+
		\frac{\nu}{2}\ltwo{x_{k+1}-x_k}^2 
		+
		\frac{\nu\stepsize ^2}{2}\ltwo{g_k}^2.
\end{align} 
By the weak convexity of $f$ and Assumption~(A1), it holds that
\begin{align}\label{proof:lem:1:inp:g:xdiff}
	\E\left [\InP{g_{k}}{x_{k-1}-x_{k}} | \mc{F}_{k-1}\right ]
	&=
	\InP{\E\left [g_{k}| \mc{F}_{k-1}\right ]}{x_{k-1}-x_{k}} 
	\nonumber\\
	&\hspace{0.0cm}
	\leq 
		f(x_{k-1}) - f(x_{k})  + \frac{\rho}{2}\ltwo{x_{k}-x_{k-1}}^2.
\end{align}
Next using  the optimality condition of $x_{k+1}$ in \eqref{alg:iters:def:a}:
\begin{align}\label{eq:x:opt:cond}
	\InP{x_{k+1} -x_k + \stepsize z_k}{x - x_{k+1}} \geq 0\quad \mbox{for all}\,\,\, x \in \mc{X},
\end{align}
we deduce (by selecting $x=x_k\in\mc{X}$) that 
\begin{align}\label{eq:x:opt:cond:x_next:new}
	\InP{ z_k}{ x_{k+1} -  x_{k}}
	\leq 
	- \stepsize^{-1}\ltwo{x_{k+1}-x_k}^2.
\end{align}
Therefore, by taking the conditional expectation in \eqref{eq:proof:lem:1:norm:2}, combining the result with \eqref{proof:lem:1:inp:g:xdiff} and \eqref{eq:x:opt:cond:x_next:new}, and rearranging terms, we arrive at
\begin{align}\label{eq:proof:lem:1:norm:bounded}
	(1-\mmt)f(x_{k})
	+	
	\E\left[		
		\frac{\nu}{2}
		\ltwo{p_{k+1}}^2|\mc{F}_{k-1}	\right]		
	&\leq
	(1-\mmt)f(x_{k-1})
	+
	\frac{\nu}{2}
	\ltwo{p_k}^2
	\nonumber\\
	&\hspace{0.45cm}
		- \left(\frac{1}{\stepsize} - \frac{\nu}{2}\right)\E\left[\ltwo{x_{k+1}-x_k}^2|\mc{F}_{k-1}	\right]	
	\nonumber\\
	&\hspace{0.45cm}
		+
		\frac{(1-\mmt)\rho}{2}\ltwo{x_{k}-x_{k-1}}^2 
		+
		\frac{\nu\stepsize ^2}{2}\E[\ltwo{g_k}^2|\mc{F}_{k-1}]	.
 \end{align}
 Next, we show by induction that $\E[\ltwo{z_k}^2|\mc{F}_{k-1}] \leq L^2\,\,\, \forall k\in\N$. Since $z_0\in\subdiff f(x_0,S_0)$, the base case follows directly from Assumption~(A2). Suppose the hypothesis holds for $i=0,\ldots,k$, we have 
\begin{align}\label{eq:znorm}
	\E[\ltwo{z_{k+1}}^2|\mc{F}_{k}]
	&= 
		\E\left[\ltwo{\mmt g_{k+1} +(1-\mmt)(x_k-x_{k+1})/{\stepsize}}^2\big|\mc{F}_{k}\right]		
	\nonumber\\
	&\mathop \leq \limits^{\mathrm{(a)}}
		\mmt\E\left[\ltwo{ g_{k+1}}^2|\mc{F}_{k}\right]
		+
		(1-\mmt)\E\left[\ltwo{(x_k-x_{k+1})/{\stepsize}}^2|\mc{F}_{k}\right]
	\nonumber\\
	&\mathop \leq \limits^{\mathrm{(b)}} 
		\mmt L^2
		+
		(1-\mmt)\E\left[\ltwo{z_k}^2|\mc{F}_{k-1}\right]	
	\mathop \leq \limits^{\mathrm{(c)}} L^2,
\end{align}
where $\mathrm{(a)}$ is true since $\ltwo{\cdot}^2$ is convex; $\mathrm{(b)}$ follows from Assumption~(A2) and the nonexpansiveness of $\proj{\cdot}$; and $\mathrm{(c)}$ follows from the induction hypothesis.
This together with the nonexpansiveness of $\proj{\cdot}$ imply that
\begin{align}\label{eq:xdiff:norm:L}
	\E[\ltwo{x_{k+1}-x_{k}}^2|\mc{F}_{k-1}] \leq \E[\ltwo{\stepsize z_k}^2|\mc{F}_{k-1} ]\leq \stepsize^2 L^2 \quad \mbox{for all}\,\,\, k\in\N.
\end{align}
Finally, using \eqref{eq:xdiff:norm:L} and  Assumption~(A2), we obtain 
\begin{align}\label{eq:proof:lem:1:norm:final}
	(1-\mmt)f(x_{k})
	+	
	\E\left[		
		\frac{\nu}{2}
		\ltwo{p_{k+1}}^2|\mc{F}_{k-1}	\right]		
	&\leq
	(1-\mmt)f(x_{k-1})
	+
	\frac{\nu}{2}
	\ltwo{p_k}^2
	\nonumber\\
	&\hspace{0.45cm}
		- \frac{1}{\stepsize}\E\left[\ltwo{x_{k+1}-x_k}^2|\mc{F}_{k-1}	\right]	
	\nonumber\\
	&\hspace{0.45cm}
		+
		\stepsize^2
		\left(\frac{\rho(1-\mmt)}{2}+\nu\right)L^2,
\end{align}
completing the proof.
\section{Proof of Lemma~\ref{lem:Vfunc:non-smooth}}\label{appendix:proof:lem:Vfunc:non-smooth}

Let $\bar{x}_{k}:= x_{k} + \frac{1-\mmt}{\mmt} \left(x_{k} - x_{k-1}\right)$ and define the \emph{virtual} iterates:
\begin{align}
	\hat{x}_{k}: 
	= 
	\prox{\lambda F}{\bar{x}_{k}} 
	= 
	\argmin_{x\in\R^n} 
	\left\{
		F(x)
		+ \frac{1}{2\lambda} \ltwo{x-\bar{x}_{k}}^2\right\}.
\end{align}
In view of Lemma~\ref{lem:moreau:env}, we have
\begin{align*}
	\grad{F_\lambda(\bar{x}_k)}=\lambda^{-1}(\bar{x}_k-\hat{x}_k),
\end{align*}
where $F_\lambda(\cdot)$ denotes the Moreau envelope of $F(x) = f(x) + \indcfunc{x}$.
By the definition of $F_{\lambda}(\bar{x}_{k+1})$, it holds that
\begin{align}\label{eq:proof:thrm:1:evn:1}
	F_{\lambda}(\bar{x}_{k+1})
	&=
		f(\hat{x}_{k+1})
		+
		\frac{1}{2\lambda}
		\ltwo{ x_{k+1} + \frac{1-\mmt}{\mmt} \left(x_{k+1} - x_{k}\right)- \hat{x}_{k+1}}^2			
	\nonumber\\
	&\leq
		f(\hat{x}_k)
		+
		\frac{1}{2\lambda}
		\ltwo{ x_{k+1} + \frac{1-\mmt}{\mmt} \left(x_{k+1} - x_{k}\right)- \hat{x}_k}^2		
	\nonumber\\
	&=
		f(\hat{x}_k)
		+
		\frac{1}{2\lambda\mmt^2}
		\ltwo{x_{k+1} -  \left((1-\mmt)x_{k} + \mmt \hat{x}_k\right)}^2.
\end{align}
Since $(1-\mmt)x_{k} + \mmt \hat{x}_k \in\mc{X}$ and $\proj{\cdot}$ is nonexpansive, we obtain
\begin{align}\label{eq:proof:thrm:1:norm:decomp:a}
	\frac{1}{2\lambda\mmt^2}
	\ltwo{x_{k+1} -  \left((1-\mmt)x_{k} + \mmt \hat{x}_k\right)}^2		
	&\leq
		\frac{1}{2\lambda\mmt^2}
		\ltwo{ x_k - \stepsize z_k - \left((1-\mmt)x_{k} + \mmt \hat{x}_k\right)}^2		
	\nonumber\\
	&=
		\frac{1}{2\lambda}
		\ltwo{ x_{k} + \frac{1-\mmt}{\mmt} \left(x_{k} - x_{k-1}\right)
		- \stepsize g_k 
		-\hat{x}_k}^2
	\nonumber\\
	&=
		\frac{1}{2\lambda}
		\ltwo{ \bar{x}_k - \hat{x}_k - \stepsize g_k}^2.			
\end{align}
Combining \eqref{eq:proof:thrm:1:evn:1} and \eqref{eq:proof:thrm:1:norm:decomp:a} yields
\begin{align}\label{eq:proof:thrm:1:evn:1:1}
	F_{\lambda}(\bar{x}_{k+1})
		\leq
		f(\hat{x}_{k+1})
		+
		\frac{1}{2\lambda}
			\ltwo{ \bar{x}_k - \hat{x}_k - \stepsize g_k}^2.
\end{align}
Using the definition of $\bar{x}_k$, we have
\begin{align}\label{eq:proof:thrm:1:norm:decomp:b}
	\ltwo{ \bar{x}_k - \hat{x}_k - \stepsize g_k}^2
	=
		\ltwo{ \bar{x}_k - \hat{x}_k}^2
		+
		2\stepsize \InP{ \hat{x}_k - x_{k} }{ g_k }
		+
		\frac{2\stepsize(1-\mmt)}{\mmt} \InP{ x_{k-1} -  x_{k}}{ g_k }
		+
		\stepsize^2\ltwo{g_k}^2.
\end{align}
It follows from \eqref{eq:proof:thrm:1:evn:1:1} and \eqref{eq:proof:thrm:1:norm:decomp:b} that
\begin{align}\label{eq:proof:thrm:1:evn:2}
	F_{\lambda}(\bar{x}_{k+1})
	&\leq
		f(\hat{x}_k)
		+
		\frac{1}{2\lambda}
		\ltwo{ \bar{x}_k - \hat{x}_k}^2
		+
		\frac{\stepsize}{\lambda} \InP{ \hat{x}_k - x_{k} }{ g_k }
		+
		\frac{\stepsize(1-\mmt)}{\lambda\mmt} \InP{ x_{k-1} -  x_{k}}{ g_k }
		+
		\frac{\stepsize^2\ltwo{g_k}^2}{2\lambda}
	\nonumber\\
	&=
		F_{\lambda}(\bar{x}_{k})
		+
		\frac{\stepsize}{\lambda} \InP{ \hat{x}_k - x_{k} }{ g_k }
		+
		\frac{\stepsize(1-\mmt)}{\lambda\mmt} \InP{ x_{k-1} -  x_{k}}{ g_k }
		+
		\frac{\stepsize^2\ltwo{g_k}^2}{2\lambda}.
\end{align}
We next bound the first inner product in \eqref{eq:proof:thrm:1:evn:2}. Since $f(\cdot)$ is $\rho$-weakly convex, it holds that
\begin{align}\label{eq:proof:thrm:1:inp:1}
	\E[\InP{\hat{x}_{k}-x_{k}}{g_{k}}|\mc{F}_{k-1}]
	&\leq	
		 f(\hat{x}_{k}) - f(x_{k}) + \frac{\rho}{2}\ltwo{\hat{x}_{k}-x_{k}}^2
	\nonumber\\
	&=
		F(\hat{x}_{k}) - F(x_{k}) + \frac{\rho}{2}\ltwo{\hat{x}_{k}-x_{k}}^2,
\end{align}
where the last step is true since $\indcfunc{x_k}=\indcfunc{\hat{x}_k}=0$. We also have that the function $x \mapsto F(x) + \frac{1}{2\lambda}\ltwo{x-\bar{x}_{k}}^2$ \, is  $(\lambda^{-1}-\rho)$-strongly convex with $\hat{x}_k$ being its minimizer, it follows from \cite[Theorem~5.25]{Bec17} that
\begin{align}\label{eq:strong:cvx}
	F(x_{k}) + \frac{1}{2\lambda} \ltwo{x_{k}-\bar{x}_{k}}^2
	-
	\left(	F(\hat{x}_{k}) + \frac{1}{2\lambda} \ltwo{\hat{x}_{k}-\bar{x}_{k}}^2\right)
	\geq
		\frac{\lambda^{-1}-\rho}{2}\ltwo{\hat{x}_{k}-x_{k}}^2.
\end{align}		
 We thus have
\begin{align}\label{eq:proof:thrm:1:inp:1:bounded}
	F(x_{k}) - F(\hat{x}_{k}) - \frac{\rho}{2}\ltwo{\hat{x}_{k}-x_{k}}^2
	&=
		F(x_{k}) + \frac{1}{2\lambda} \ltwo{x_{k}-\bar{x}_{k}}^2
		-
		\left(	F(\hat{x}_{k}) + \frac{1}{2\lambda} \ltwo{\hat{x}_{k}-\bar{x}_{k}}^2\right)
	\nonumber\\
	&\hspace{0.45cm}
		-
		\frac{1}{2\lambda} \ltwo{x_{k}-\bar{x}_{k}}^2
		+
		\frac{1}{2\lambda} \ltwo{\hat{x}_{k}-\bar{x}_{k}}^2
		- 
		\frac{\rho}{2}\ltwo{\hat{x}_{k}-x_{k}}^2
	\nonumber\\
	&\hspace{0.0cm}
	\mathop \geq \limits^{\mathrm{(a)}}
		\frac{\lambda^{-1}-\rho}{2}\ltwo{\hat{x}_{k}-x_{k}}^2
		-
		\frac{1}{2\lambda} \ltwo{x_{k}-\bar{x}_{k}}^2
	\nonumber\\
	&\hspace{0.450cm}
		+
		\frac{1}{2\lambda} \ltwo{\hat{x}_{k}-\bar{x}_{k}}^2
		- 
		\frac{\rho}{2}\ltwo{\hat{x}_{k}-x_{k}}^2
	\nonumber\\
	&\hspace{.0cm}
	\mathop = \limits^{\mathrm{(b)}}
		\frac{\lambda}{2} \ltwo{\grad{F_{\lambda}(\bar{x}_{k})}}^2
		+
		\frac{\lambda^{-1}-2\rho}{2}\ltwo{\hat{x}_{k}-x_{k}}^2
		-
		\frac{1}{2\lambda} \ltwo{x_{k}-\bar{x}_{k}}^2
	\nonumber\\
	&\hspace{.0cm}
	\mathop \geq \limits^{\mathrm{(c)}}
		\frac{\lambda}{2} \ltwo{\grad{F_{\lambda}(\bar{x}_{k})}}^2
		-
		\frac{1}{2\lambda} \ltwo{x_{k}-\bar{x}_{k}}^2,
\end{align}
where $\mathrm{(a)}$ is due to \eqref{eq:strong:cvx},  $\mathrm{(b)}$ follows from the definition of $\grad{F_{\lambda}(\bar{x}_{k})}$, and $\mathrm{(c)}$ holds since $\lambda^{-1}\geq 2\rho$.
 For the second inner product in \eqref{eq:proof:thrm:1:evn:2}, we have
\begin{align}\label{eq:proof:thrm:1:inp:2}
	\E[\InP{ x_{k-1} -  x_{k}}{ g_k }| \mc{F}_{k-1}]
	\leq 
		f(x_{k-1}) - f(x_k)
		+
		\frac{\rho}{2}\ltwo{x_{k-1} -  x_{k}}^2.
\end{align}

\noindent Therefore, combining \eqref{eq:proof:thrm:1:inp:1}, \eqref{eq:proof:thrm:1:inp:1:bounded},  \eqref{eq:proof:thrm:1:inp:2} and plugging the result into \eqref{eq:proof:thrm:1:evn:2} yield
\begin{align*}
	 \E[F_{\lambda}(\bar{x}_{k+1})| \mc{F}_{k-1}] + \frac{\stepsize(1-\mmt)}{\lambda\mmt}  f(x_k)
	&\leq
		 F_{\lambda}(\bar{x}_{k}) 
		+
		\frac{\stepsize(1-\mmt)}{\lambda\mmt}  f(x_{k-1})
		-
		\frac{\stepsize}{2}		
		\ltwo{\grad{F_{\lambda}(\bar{x}_k)}}^2
	\nonumber\\
	&\hspace{-2cm}
		+
		\frac{\stepsize}{2\lambda^2} \ltwo{x_{k}-\bar{x}_{k}}^2
		+
		\frac{\rho\stepsize(1-\mmt)}{2\lambda\mmt}
		\ltwo{x_{k-1} -  x_{k}}^2
		+
		\frac{\stepsize^2\E[\ltwo{g_k}^2| \mc{F}_{k-1}]}{2\lambda}.
\end{align*}
Now, using \eqref{eq:xdiff:norm:L}, Assumption~(A2), the definition of $\bar{x}_k$, and the fact that $\mmt=\nu\stepsize$ give
\begin{align}\label{eq:proof:thrm:1:evn:3}
	 \E[F_{\lambda}(\bar{x}_{k+1})| \mc{F}_{k-1}] + \frac{1-\mmt}{\lambda\nu}  f(x_k)
	&\leq
		 F_{\lambda}(\bar{x}_{k}) 
		+
		\frac{1-\mmt}{\lambda\nu}  f(x_{k-1})
		-
		\frac{\stepsize}{2}		
		\ltwo{\grad{F_{\lambda}(\bar{x}_k)}}^2
	\nonumber\\
	&\hspace{-2cm}
		+
		\frac{(1-\mmt)^2}{2\lambda^2\nu^2}\frac{1}{\stepsize} \ltwo{x_{k}-x_{k-1}}^2
		+
		\stepsize^2
		\left(\frac{\rho(1-\mmt)}{2\nu}+1\right)
		\frac{L^2}{2\lambda}
\end{align}
For simplicity, define $\xi=(1-\mmt)/\nu$. Finally, to form a telescoping sum, we simply multiply both sides of eq.~\eqref{eq:lem:xdiff:non-smooth} in Lemma~\eqref{lem:xdiff:non-smooth} by $\xi^2/(2\lambda^2)$ and combine the result with \eqref{eq:proof:thrm:1:evn:3} to get
\begin{align*}
	 &\E\left[F_{\lambda}(\bar{x}_{k+1})
	 +
 	 \frac{\nu\xi^2}{4\lambda^2}
 	 \ltwo{p_{k+1}}^2
 	 + 
 	 \frac{\xi^2}{2\stepsize\lambda^2}
	 \ltwo{x_{k+1}-x_k}^2 \Big|\mc{F}_{k-1}	\right]	
	 +
	 \left(\frac{(1-\mmt)\xi^2}{2\lambda^2}+\frac{\xi}{\lambda}\right)  f(x_k)
	\nonumber\\
	&\hspace{1cm}
	\leq
		 F_{\lambda}(\bar{x}_{k}) 
		+
		\frac{\nu\xi^2}{4\lambda^2}
	 	 \ltwo{p_{k}}^2
	 	+
	 	\frac{\xi^2}{2\stepsize\lambda^2}\ltwo{x_{k}-x_{k-1}}^2
		+
		\left(\frac{(1-\mmt)\xi^2}{2\lambda^2}+\frac{\xi}{\lambda}\right)  f(x_{k-1})
	\nonumber\\
	&\hspace{1.45cm}
		-
		\frac{\stepsize}{2}		
		\ltwo{\grad{F_{\lambda}(\bar{x}_k)}}^2	
		+
		\frac{\gamma\stepsize^2 L^2}{2\lambda},
\end{align*}
where 
\begin{align*}
	\gamma 
	= 
	\frac{\xi^2\left(\frac{\rho(1-\mmt)}{2}+\nu\right)}{\lambda} 
	+
	\frac{\rho\xi}{2}+1.
\end{align*}
The proof is complete.

\section{Proof of Lemma~\ref{lem:Wfunc:smooth}}
\label{appendix:proof:Wfunc:smooth}

In this section, we prove Lemma~\ref{lem:Wfunc:smooth}. We begin with the following lemma.
\begin{lemma}\label{lem:xdiff:smooth}
Let  Assumptions~(A1) and (A3) hold. Let $\stepsize\in(0, 1/\rho)$ and $\mmt= \nu\stepsize$ for some constant $\nu>0$ such that $\mmt\in (0,1]$. Then, for any $k\in\N$, the iterates generated by procedure \eqref{alg:iters:def} satisfy:
\begin{align*}
	\E\left[f(x_{k+1}) + \frac{\varphi_{k+1}}{\nu\stepsize^2} \Big| \mc{F}_{k}\right ]
	&\leq
		 f(x_k)+ \frac{\varphi_{k}}{\nu\stepsize^2}
		-(\stepsize -\frac{\rho\stepsize^2 }{2})
		\ltwo{d_{k+1}}^2
	 	+
	 	\frac{1}{2\nu}\E\left [ \ltwo{ z_{k} -  z_{k+1}}^2| \mc{F}_{k}\right ].
\end{align*}
\begin{proof}

We first rewrite the update of $x_{k+1}$ in \eqref{alg:iters:def:a} on the following form:
\begin{align}\label{alg:iters:def:a:rewrite}
	x_{k+1} = \argmin_{x\in\mc{X}} \left\{ \InP{\stepsize z_k-x_k}{x} + \frac{1}{2} \ltwo{x}^2\right\}.
\end{align}
Let  $h(x)=\frac{1}{2}\ltwo{x}^2+\indcfunc{x}$ and define its convex conjugate 
$h^*(y) = \sup_{x\in\mc{X}}\{\InP{y}{x} - h(x)\}$. It is well-known that $\grad{h^*}$ is $1$-Lipschitz with  gradient $\grad{h^*(y)} = \argmax_{x\in\mc{X}}\{\InP{y}{x} - h(x)\}$ \cite[Chapter~X]{HL93}.  
Therefore, the update formula \eqref{alg:iters:def:a:rewrite} implies that $\grad{h^*(x_{k} - \stepsize z_{k})} = x_{k+1}$. By the smoothness of $h^*$, we have
\begin{align*}
	 h^*(x_{k+1} - \stepsize z_{k+1})  
	 &\leq 
	 h^*(x_{k} - \stepsize z_{k}) 
	 	+
	 	\InP{x_{k+1}}{ x_{k+1} - \stepsize z_{k+1} - x_{k} + \stepsize z_{k}}
 	\nonumber\\
	&\hspace{0.45cm}
		+
		\frac{1}{2}\ltwo{x_{k} - \stepsize z_{k} - x_{k+1} + \stepsize z_{k+1}}^2.
\end{align*}
Let $\varphi_k:=h^*(x_{k} - \stepsize z_{k}) - \half \ltwo{x_{k}}^2 + \stepsize \InP{x_{k}}{ z_{k}}$, it then follows that  
\begin{align}\label{eq:varphi:diff}
	\varphi_{k+1} - \varphi_k 
	&=
		h^*(x_{k+1} - \stepsize z_{k+1})  
		-
		\half \ltwo{x_{k+1}}^2
		+
		\stepsize\InP{x_{k+1}}{  z_{k+1}}
	\nonumber\\
	&\hspace{0.45cm}
		-
		h^*(x_{k} - \stepsize z_{k}) 
		+
		\half \ltwo{x_{k}}^2
		-
		\stepsize\InP{x_{k}}{  z_{k}}				
	\nonumber\\
	&\leq
	 	\InP{x_{k+1}}{ x_{k+1} - \stepsize z_{k+1} - x_{k} + \stepsize z_{k}}
		+
		\frac{1}{2}\ltwo{x_{k} - \stepsize z_{k} - x_{k+1} + \stepsize z_{k+1}}^2
	\nonumber\\
	&\hspace{0.45cm}			
		-
		\half \ltwo{x_{k+1}}^2
		+
		\stepsize \InP{x_{k+1}}{ z_{k+1}}
		+
		\half \ltwo{x_{k}}^2
		-
		\stepsize\InP{x_{k}}{  z_{k}}.
\end{align}
We next apply the identity
\begin{align}\label{eq:norm:identity}
	\InP{a}{b} = \half \ltwo{a}^2 + \half \ltwo{b}^2 - \half \ltwo{a-b}^2
\end{align}
to obtain
\begin{align}\label{proof:lema:1:inp:1:a}
	\InP{x_{k}}{\stepsize  z_{k}}
	&=
	\half
	\ltwo{x_{k}}^2
	+
	\frac{\stepsize ^2}{2}
	\ltwo{z_{k}}^2
	-
	\half
	\ltwo{x_{k}-\stepsize z_{k}}^2
	\\ \label{proof:lema:1:inp:1:b}
	\InP{x_{k+1}}{\stepsize  z_{k+1}}
	&=
	\half
	\ltwo{x_{k+1}}^2
	+
	\frac{\stepsize ^2}{2}
	\ltwo{z_{k+1}}^2
	-
	\half
	\ltwo{x_{k+1}-\stepsize z_{k+1}}^2.
\end{align}
Using the identity \eqref{eq:norm:identity} again, we get
\begin{align}\label{proof:lema:1:inp:2}
	\half
	\ltwo{x_{k}-\stepsize z_{k}}^2
	-
	\half
	\ltwo{x_{k+1}-\stepsize z_{k+1}}^2
	&=
		\InP{x_k - \stepsize  z_k}{x_{k} - \stepsize z_{k} - x_{k+1} + \stepsize z_{k+1}}
	\nonumber\\
	&\hspace{0.45cm}
		-
		\frac{1}{2}\ltwo{x_{k} - \stepsize z_{k} - x_{k+1} + \stepsize z_{k+1}}^2.
\end{align}
Plugging \eqref{proof:lema:1:inp:1:a}--\eqref{proof:lema:1:inp:2} into \eqref{eq:varphi:diff} yields
\begin{align}\label{eq:varphi:diff:2}
	\varphi_{k+1} - \varphi_k 
	&\leq
		\frac{\stepsize ^2}{2}
			\ltwo{z_{k+1}}^2
		-
		\frac{\stepsize ^2}{2}
			\ltwo{z_{k}}^2
		+
	 	\InP{x_{k+1}-x_k + \stepsize  z_k}{ x_{k+1} - x_{k} + \stepsize z_{k} - \stepsize z_{k+1}}
	\nonumber\\
	&=
		\frac{\stepsize ^2}{2}
			\ltwo{z_{k+1}}^2
		-
		\frac{\stepsize ^2}{2}
			\ltwo{z_{k}}^2
		+
		\stepsize^2\InP{ z_k}{ z_{k} -  z_{k+1}}	
	\nonumber\\
	&\hspace{0.45cm}
		+
	 	\InP{x_{k+1}-x_k + \stepsize  z_k}{ x_{k+1} - x_{k}}		
	 	+
	 	\stepsize\InP{x_{k+1}-x_k }{ z_{k} -  z_{k+1}}	
\end{align}
By identity \eqref{eq:norm:identity}, we have
\begin{align}\label{proof:lema:1:inp:3:c}
	\frac{1}{2}
			\ltwo{z_{k+1}}^2
		-
		\frac{1}{2}
			\ltwo{z_{k}}^2
	+
	\InP{ z_k}{ z_{k} - z_{k+1}}
	=
		\frac{1}{2} \ltwo{z_{k} - z_{k+1}}^2.
\end{align}
Note that the optimality condition of $x_{k+1}$ in \eqref{eq:x:opt:cond}
implies that
\begin{align}\label{eq:x:opt:cond:x_next}
	\InP{x_{k+1}-x_k + \stepsize  z_k}{ x_{k+1} -  x_{k} }
	\leq 0.
\end{align}
By the definition of $z_{k+1}$, it holds that
\begin{align}\label{proof:lema:1:inp:3:a}
	\InP{x_{k+1}-x_k}{ \stepsize z_{k} - \stepsize z_{k+1}}
	&=
		\InP{x_{k+1}-x_k}{ \stepsize z_{k} }
		-
		 \stepsize \mmt 
		\InP{x_{k+1}-x_k}{ g_{k+1} }
	\nonumber\\
	&\hspace{0.45cm}
		+
		(1-\mmt )
		\ltwo{x_{k+1}-x_k}^2.
\end{align}
We can also deduce from \eqref{eq:x:opt:cond:x_next} that $\InP{x_{k+1}-x_k}{ \stepsize z_{k} }\leq -\ltwo{x_{k+1}-x_k}^2$, and hence \eqref{proof:lema:1:inp:3:a} can be further bounded by
\begin{align}\label{proof:lema:1:inp:3:b}
	\InP{x_{k+1}-x_k}{ \stepsize z_{k} - \stepsize z_{k+1}}
	&\leq
		-\mmt 
		\ltwo{x_{k+1}-x_k}^2
		+
		\stepsize \mmt 
		\InP{x_k-x_{k+1}}{ g_{k+1} }.
\end{align}
Thus, by combining \eqref{eq:varphi:diff:2}, \eqref{proof:lema:1:inp:3:c},  \eqref{eq:x:opt:cond:x_next}, and \eqref{proof:lema:1:inp:3:b}, we arrive at
\begin{align}\label{eq:varphi:diff:3}
	\varphi_{k+1} 
	&\leq
		 \varphi_k 
		-\mmt 
		\ltwo{x_{k+1}-x_k}^2
		+
	 	\stepsize \mmt  			
	 	\InP{x_k-x_{k+1}}{ g_{k+1} }
	 	+
	 	\frac{\stepsize^2}{2} \ltwo{ z_{k} -  z_{k+1}}^2.
\end{align}
Next, by the weak convexity of $f$ and Assumption~(A1), we have
\begin{align}\label{proof:lema:1:inp:g:xdiff}
	f(x_{k+1})
	\leq 
		f(x_k) - \E\left [\InP{x_k-x_{k+1}}{g_{k+1}} | \mc{F}_{k}\right ]  + \frac{\rho}{2}\ltwo{x_{k+1}-x_k}^2.
\end{align}
Thus, multiplying both sides of \eqref{eq:varphi:diff:3} by $1/(\mmt\stepsize)=1/(\nu\stepsize^2)$, taking the conditional expectation, and adding the result to \eqref{proof:lema:1:inp:g:xdiff} give
\begin{align}\label{eq:varphi:diff:4}
	\E\left[f(x_{k+1}) + \frac{\varphi_{k+1}}{\nu\stepsize^2} \Big| \mc{F}_{k}\right ]
	&\leq
		 f(x_k)+ \frac{\varphi_{k}}{\nu\stepsize^2}
		-(\frac{1}{\stepsize} -\frac{\rho }{2})
		\ltwo{x_{k+1}-x_k}^2
	\nonumber\\
	&\hspace{0.5cm}
	 	+
	 	\frac{1}{2\nu}\E\left [ \ltwo{ z_{k} -  z_{k+1}}^2| \mc{F}_{k}\right ].
\end{align}
Using the definition of $d_{k+1}$ completes the proof.

\end{proof}
\end{lemma}

The next lemma bounds the term $\frac{1}{2\nu}\E\left [ \ltwo{ z_{k} -  z_{k+1}}^2| \mc{F}_{k}\right ]$ in \eqref{eq:varphi:diff:4}. 
\begin{lemma} \label{lem:z:diff}
Let $\xi=(1-\mmt)/\nu$. Under the same setting of Lemma~\ref{lem:xdiff:smooth}, we have
\begin{align*}
	\E\left[\frac{1}{2\nu}\ltwo{z_{k+1}-z_k}^2|\mc{F}_k\right]
	&\leq
		f(x_k)-f(x_{k+1})
		+		
		\frac{\xi}{2}\ltwo{d_k}^2	- \frac{\xi}{2} \ltwo{d_{k+1}}^2	
	\nonumber\\	 	
	&\hspace{0.45cm}
		-\left(\stepsize-\frac{\stepsize^3\rho^2+3\rho\stepsize^2}{2}\right)\ltwo{d_{k+1}}^2
		+
		4\nu\stepsize^2\sigma^2.
\end{align*}

\begin{proof}
Let $\Delta_{k}:=g_k-\grad{f(x_k)}$, we have
\begin{align}\label{proof:smooth:z:diff:1}
	\ltwo{z_{k+1}-z_k}^2
	&=
		\ltwo{
			\mmt \left(g_{k+1}-g_k\right) 
			+
			\left(1-\mmt\right)\left(d_{k+1}-d_k\right)
		}^2
	\nonumber\\
	&\hspace{0cm}
	=
		\ltwo{
			\mmt \left(\Delta_{k+1}-\Delta_k\right) 
			+
			\mmt \left(\grad{f(x_{k+1})}-\grad{f(x_{k})}\right) 
			+
			\left(1-\mmt\right)\left(d_{k+1}-d_k\right)
		}^2
	\nonumber\\
	&\hspace{0cm}	
	=
	\mmt^2\ltwo{\Delta_{k+1}-\Delta_k}^2
	+
	\ltwo{
		\mmt \left(\grad{f(x_{k+1})}-\grad{f(x_{k})}\right) 
		+
		\left(1-\mmt\right)\left(d_{k+1}-d_k\right)}^2
	\nonumber\\
	&\hspace{.45cm}	
	+
		2\mmt
		\InP{
			\Delta_{k+1}-\Delta_k
		}{
			\mmt \left(\grad{f(x_{k+1})}-\grad{f(x_{k})}\right) 
			+
			\left(1-\mmt\right)\left(d_{k+1}-d_k\right)
		}.	
\end{align}
First, by Assumption (A3), we have
\begin{align}\label{proof:smooth:z:diff:Delta}
	\mmt^2\E\left[ \ltwo{\Delta_{k+1}-\Delta_k}^2 | \mc{F}_{k} \right]
	\leq
		\mmt^2\E\left[ 2\ltwo{\Delta_{k+1}}^2 + 2\ltwo{\Delta_k}^2 | \mc{F}_{k} \right]
	\leq 4\mmt^2\sigma^2.
\end{align}
Define $p_k:= \mmt \grad{f(x_{k})}	+ \left(1-\mmt\right)d_k$ and let $T$ be the last term in \eqref{proof:smooth:z:diff:1}, then
\begin{align}\label{proof:smooth:z:diff:T}
	T=
	2\mmt
	\InP{
		\Delta_{k+1}-\Delta_k
	}{
		p_{k+1}-p_k
	}.	
\end{align}
Since $p_{k+1}-p_k$ is a function of $\mc{F}_k = \sigma(\statrv_0, \ldots, \statrv_k)$, it follows from Assumption~(A1) that
\begin{align}\label{proof:smooth:cond:expect:a}
	\E\left[\InP{\Delta_{k+1}}{p_{k+1}-p_k}|\mc{F}_{k}\right]
	&=
		\E\left[\InP{g_{k+1}-\grad{f(x_{k+1})}}{p_{k+1}-p_k}|\mc{F}_{k}\right]				
	\nonumber\\	
	&=
		\InP{\E\left[g_{k+1}|\mc{F}_{k}\right]-\grad{f(x_{k+1})}}{p_{k+1}-p_k}		
	= 
		0.
\end{align}
Similarly, conditioned on $\mc{F}_k$, both terms in the inner product of $\E\left[\InP{\Delta_{k}}{p_{k+1}-p_k}|\mc{F}_{k}\right]$ are deterministic, and hence
\begin{align}\label{proof:smooth:cond:expect:b}
	-\E\left[\InP{\Delta_{k}}{p_{k+1}-p_k}|\mc{F}_{k}\right]
	=
		-\E\left[\InP{\Delta_{k}}{p_{k+1}-p_k}|\mc{F}_{k-1}\right].
\end{align} 
Similar to \eqref{proof:smooth:cond:expect:a}, we obtain
\begin{align}\label{proof:smooth:cond:expect:b:1}
	\E\left[\InP{\Delta_{k}}{p_k}|\mc{F}_{k-1}\right]
	=
		\InP{\E\left[g_{k}|\mc{F}_{k-1}\right]-\grad{f(x_{k})}}{p_k}		
	= 
		0.
\end{align}
To bound the remaining term $\E\left[\InP{\Delta_{k}}{p_{k+1}}|\mc{F}_{k-1}\right]$, we  introduce the following \emph{virtual} iterate:
\begin{align*}
	x'_{k+1} = \proj{x_k - \stepsize p_k },
\end{align*}
and define 
\begin{align*}
	p'_{k+1} = \mmt\grad{f(x'_{k+1})} + (1-\mmt)\frac{x_k-x'_{k+1}}{\stepsize}.
\end{align*}
By definitions, both $x'_{k+1}$ and $p'_{k+1}$ only depend on $\mc{F}_{k-1}$, it follows that
\begin{align*}
	-\E\left[\InP{\Delta_{k}}{p_{k+1}}|\mc{F}_{k-1}\right]
	&=
		-\E\left[\InP{\Delta_{k}}{p'_{k+1}}|\mc{F}_{k-1}\right]
		+
		\E\left[\InP{\Delta_{k}}{p'_{k+1}-p_{k+1}}|\mc{F}_{k-1}\right]
	\nonumber\\
	&=
		0 + \E\left[\InP{\Delta_{k}}{p'_{k+1}-p_{k+1}}|\mc{F}_{k-1}\right]
	\nonumber\\
	&\leq
		\E\left[\ltwo{\Delta_{k}}\ltwo{p_{k+1}-p'_{k+1}}|\mc{F}_{k-1}\right].		
\end{align*}
We have 
\begin{align}\label{proof:smooth:p:normdiff}
	\ltwo{p_{k+1}-p'_{k+1}} 
	&=
		\ltwo{\mmt(\grad{f(x_{k+1})}-\grad{f(x'_{k+1})}) + (1-\mmt)\frac{x'_{k+1}-x_{k+1}}{\stepsize}}
	\nonumber\\
	&\mathop \leq \limits^{\mathrm{(a)}}
		\mmt\ltwo{(\grad{f(x_{k+1})}-\grad{f(x'_{k+1})})}
		+
		\frac{1-\mmt}{\stepsize}\ltwo{x'_{k+1}-x_{k+1}}
	\nonumber\\
	&\mathop \leq \limits^{\mathrm{(b)}}
		\left(\mmt\rho + \frac{1-\mmt}{\stepsize}\right)\ltwo{x_{k+1}-x'_{k+1}}
	\nonumber\\
	&\mathop \leq \limits^{\mathrm{(c)}}
		\stepsize\mmt\left(\mmt\rho + \frac{1-\mmt}{\stepsize}\right)\ltwo{g_k-\grad{f(x_k)}}
	\nonumber\\
	&\mathop \leq \limits^{\mathrm{(d)}}
		\mmt\ltwo{\Delta_k},
\end{align}
where $\mathrm{(a)}$ holds since $\ltwo{\cdot}$ is convex; $\mathrm{(b)}$ follows since $\grad{f}$ is $\rho$-Lipschitz; $\mathrm{(c)}$ follows from the definition of $x'_{k+1}$, $p_k$, and the nonexpansiveness of $\proj{\cdot}$; and $\mathrm{(d)}$ is true since $\stepsize\in (0, 1/\rho)$. We thus have 
\begin{align}\label{proof:smooth:cond:expect:b:3}
	-\E\left[\InP{\Delta_{k}}{p_{k+1}}|\mc{F}_{k-1}\right]
	&\leq
		\mmt\E\left[\ltwo{\Delta_k}^2|\mc{F}_{k-1}\right]
	\leq
	\mmt\sigma^2.
\end{align}
Therefore, taking the conditional expectation in \eqref{proof:smooth:z:diff:T} and using \eqref{proof:smooth:cond:expect:a}, \eqref{proof:smooth:cond:expect:b}, \eqref{proof:smooth:cond:expect:b:1}, \eqref{proof:smooth:cond:expect:b:3} yield
\begin{align}\label{proof:smooth:cond:expect:b:final}
	\E\left[T|\mc{F}_{k}\right]
	\leq
		2\mmt^2\sigma^2.
\end{align}

Plugging \eqref{proof:smooth:z:diff:Delta} and \eqref{proof:smooth:cond:expect:b:final} into \eqref{proof:smooth:z:diff:1}, we arrive at
\begin{align}\label{proof:smooth:z:diff:2}
	\E\left[\ltwo{z_{k+1}-z_k}^2|\mc{F}_k\right]
	&\leq
		\ltwo{
			\mmt \left(\grad{f(x_{k+1})}-\grad{f(x_{k})}\right) 
			+
			\left(1-\mmt\right)\left(d_{k+1}-d_k\right)}^2		
		+
		6\mmt^2\sigma^2
	\nonumber\\
	&\leq	
		\mmt\ltwo{\grad{f(x_{k+1})}-\grad{f(x_{k})} }^2
		+
		\left(1-\mmt\right)
		\ltwo{d_{k+1}-d_k}^2
		+ 
		6\mmt^2\sigma^2
	\nonumber\\
	&\leq
		\mmt\rho^2\ltwo{ x_{k+1}-x_{k}}^2
		+
		\left(1-\mmt\right)
		\ltwo{d_{k+1}-d_k}^2
		+
		6\mmt^2\sigma^2.
	\nonumber\\
	&=
		\mmt\rho^2\stepsize^2\ltwo{d_{k+1}}^2
		+
		\left(1-\mmt\right)
		\ltwo{d_{k+1}-d_k}^2
		+
		6\mmt^2\sigma^2.	
\end{align} 
where we used the convexity of $\ltwo{\cdot}^2$ in the second inequality.
We now decompose and bound the term $\ltwo{d_{k+1}-d_k}^2$ as follows
\begin{align}\label{proof:smooth:d:diff:decom}
	\ltwo{d_{k+1}-d_k}^2
	=
	\ltwo{d_k}^2
	-
	\ltwo{d_{k+1}}^2
	+
	2\InP{d_{k+1}}{d_{k+1}-d_k}.
\end{align}
First, it follows from the optimality condition of $x_{k+1}$ in \eqref{eq:x:opt:cond} and the definition of $d_k$ that
\begin{align*}
	\InP{-d_{k+1} + \mmt g_k + (1-\mmt) d_k}{-d_{k+1}}\leq 0.
\end{align*}
Thus, 
\begin{align}\label{proof:smooth:d:diff:decom:1}
	\InP{d_{k+1}}{d_{k+1}-d_k}
	&\leq
		\ltwo{d_{k+1}}^2
		+
		\frac{1}{1-\mmt}
		\InP{-d_{k+1} + \mmt g_k}{d_{k+1}}
	\nonumber\\	
	&=
		-\frac{\mmt}{1-\mmt}\ltwo{d_{k+1}}^2
		+
		\frac{\mmt}{\stepsize(1-\mmt)}
		\InP{g_k}{x_k-x_{k+1}}.
\end{align}
We have 
\begin{align*}
	\InP{g_k}{x_k-x_{k+1}}
	&=
		\InP{g_k-\grad{f(x_k)}}{x_k-x_{k+1}}
	\nonumber\\
	&\hspace{0.45cm}+
		\InP{\grad{f(x_{k+1})}}{x_k-x_{k+1}}
	\nonumber\\
	&\hspace{0.45cm}+
		\InP{\grad{f(x_k)}-\grad{f(x_{k+1})}}{x_k-x_{k+1}}.
\end{align*}
For the first term, since
\begin{align*}
		\E\left[\InP{g_k-\grad{f(x_k)}}{x_k-x_{k+1}}|\mc{F}_{k}\right]		
		&=
			\E\left[\InP{g_k-\grad{f(x_k)}}{x_k-x'_{k+1}}|\mc{F}_{k-1}\right]	
		\nonumber\\
		&\hspace{0.45cm}
			+
			\E\left[\InP{g_k-\grad{f(x_k)}}{x'_{k+1}-x_{k+1}}|\mc{F}_{k-1}\right]	
		\nonumber\\
		&=
			\E\left[\InP{g_k-\grad{f(x_k)}}{x'_{k+1}-x_{k+1}}|\mc{F}_{k-1}\right],
\end{align*}
following the same steps leading to \eqref{proof:smooth:p:normdiff}, it can be bounded as
\begin{align}\label{proof:smooth:g:inp}
	\E\left[\InP{g_k-\grad{f(x_k)}}{x_k-x_{k+1}}|\mc{F}_{k}\right]
	 \leq
		\stepsize\mmt\E\left[\ltwo{\Delta_{k}}^2|\mc{F}_{k-1}\right]		
		\leq \stepsize\mmt \sigma^2.
\end{align}
By the smoothness of $\grad{f}$, we also have 
\begin{align}\label{proof:smooth:g:inp:2a}
	\InP{\grad{f(x_{k+1})}}{x_k-x_{k+1}}
	&\leq 
		f(x_k) - f(x_{k+1})  + \frac{\rho}{2}\ltwo{x_{k+1}-x_k}^2
	\\
	\label{proof:smooth:g:inp:2b}
	\InP{\grad{f(x_k)}-\grad{f(x_{k+1})}}{x_k-x_{k+1}}
	&\leq
		\rho\ltwo{x_{k+1}-x_k}^2.
\end{align}
Hence, it follows from \eqref{proof:smooth:d:diff:decom:1}--\eqref{proof:smooth:g:inp:2b} that
\begin{align}\label{proof:smooth:d:diff:decom:2}
	\InP{d_{k+1}}{d_{k+1}-d_k}
	&\leq
		-\frac{\mmt}{1-\mmt}\ltwo{d_{k+1}}^2
	\nonumber\\	
	&+
		\frac{\mmt}{\stepsize(1-\mmt)}
		\left[
			f(x_k) - f(x_{k+1})
			+
			\frac{3\rho\stepsize^2}{2}\ltwo{d_{k+1}}^2
			+
			\stepsize\mmt \sigma^2
		\right],
\end{align}
where we also used the fact that $\ltwo{x_{k+1}-x_k}^2=\stepsize^2\ltwo{d_{k+1}}^2$.

By multiplying both sides of \eqref{proof:smooth:d:diff:decom} by $1-\mmt$ and combining with \eqref{proof:smooth:d:diff:decom:2} give
\begin{align}\label{proof:smooth:d:diff:decom:final}
	(1-\mmt)\ltwo{d_{k+1}-d_k}^2
	&=
		(1-\mmt)\left(\ltwo{d_k}^2-\ltwo{d_{k+1}}^2\right)	
	\nonumber\\	
	&\hspace{0.45cm}
		-
		\left({2\mmt} - 3\rho\stepsize\mmt\right)\ltwo{d_{k+1}}^2
	\nonumber\\	
	&\hspace{0.45cm}			
		+
		\frac{2\mmt}{\stepsize}
		\left(f(x_k) - f(x_{k+1})\right)
		+
		2\mmt^2 \sigma^2.
\end{align}
Plugging \eqref{proof:smooth:d:diff:decom:final} into \eqref{proof:smooth:z:diff:2} yields
\begin{align}\label{proof:smooth:z:diff:3}
	\E\left[\ltwo{z_{k+1}-z_k}^2|\mc{F}_k\right]
	&\leq
		(1-\mmt)\left(\ltwo{d_k}^2-\ltwo{d_{k+1}}^2\right)	
		+
		\frac{2\mmt}{\stepsize}\left(f(x_k)-f(x_{k+1})\right)
	\nonumber\\
	&\hspace{0.45cm}		
		-
		\left({2\mmt} - 3\rho\stepsize\mmt-\mmt\rho^2\stepsize^2\right)\ltwo{d_{k+1}}^2
		+
		8\mmt^2\sigma^2.
\end{align} 
Multiplying both sides of \eqref{proof:smooth:z:diff:3} by $1/(2\nu)=\stepsize/(2\mmt)$ and noting that $\xi=(1-\mmt)/\nu$ complete the proof.
\end{proof}
\end{lemma}

Having established Lemmas~\ref{lem:xdiff:smooth} and \ref{lem:z:diff}, the result of Lemma~\ref{lem:Wfunc:smooth} follows immediately from the definition of the function $W_k$ and the fact that $\stepsize\in\left(0,1/(4\rho)\right]$.

\section{Proof of Theorem~\ref{thrm:smooth:complexity}}\label{appendix:proof:thrm:smooth:complexity}

We first show by induction that $\E[\ltwo{z_k}^2]\leq \sigma^2+G^2$ for all $k\in\N$. Since $z_0=\grad{f(x_0,S_0)}$ and $\ltwo{\grad{f(x_0)}}\leq G$, it follows from Assumptions (A1) and (A3) that
\begin{align*}
	\E[\ltwo{z_0}^2 ]
	&= 
	\E[\ltwo{\grad{f(x_0,S_0)}-\grad{f(x_0)}+\grad{f(x_0)}}^2]
	\nonumber\\
	&= 
		\E[\ltwo{\grad{f(x_0,S_0)}-\grad{f(x_0)}}^2]
		+
		\ltwo{\grad{f(x_0)}}^2
	\nonumber\\
	&\leq 
		\sigma^2 + G^2.
\end{align*}
Suppose that $\E[\ltwo{z_i}^2] \leq \sigma^2+G^2$ for $i=0,\ldots,k$, we have
\begin{align*}
	\E[\ltwo{z_{k+1}}^2 ]
	&= 
	\E\left[\ltwo{\mmt\left(g_{k+1}-\grad{f(x_{k+1})}\right) + \mmt\grad{f(x_{k+1})} + (1-\mmt)\frac{x_k-x_{k+1}}{\stepsize}}^2\right]
	\nonumber\\
	&= 
		\mmt^2\E[\ltwo{g_{k+1}-\grad{f(x_{k+1})}}^2]
		+
		\E\left[\ltwo{ \mmt\grad{f(x_{k+1})} + (1-\mmt)\frac{x_k-x_{k+1}}{\stepsize}}^2\right]
	\nonumber\\
	&\leq 
		\mmt^2\sigma^2 
		+ \mmt\ltwo{ \grad{f(x_{k+1})}}^2
		+
		(1-\mmt)
		\E\left[\ltwo{ \frac{x_k-x_{k+1}}{\stepsize}}^2\right],
\end{align*}
where we used the convexity of $\ltwo{\cdot}^2$ in the last step. 
Since $x_{k+1}=\proj{x_k-\stepsize z_k}$, the nonexpansiveness of $\proj{\cdot}$ and the induction hypothesis imply that
\begin{align*}
	\E\left[\ltwo{ \frac{x_k-x_{k+1}}{\stepsize}}^2\right]
	\leq 
	\E\left[	\ltwo{z_k}^2 \right]
	\leq \sigma^2 + G^2.
\end{align*}
Since $\mmt\in(0,1]$, we have $\mmt^2\leq \mmt$, and hence
$\E[\ltwo{z_{k+1}}^2 ] \leq \sigma^2+G^2$, as desired. From this, together with the fact that $\E[\ltwo{g_k}^2]\leq \sigma^2+G^2$, the theorem follows immediately from replacing $L^2$ by $\sigma^2+G^2$ everywhere in the proofs of Lemmas~\ref{lem:xdiff:non-smooth}--\ref{lem:Vfunc:non-smooth} and Theorem~\ref{thrm:non-smooth:complexity}.

\section{Proof of Theorem~\ref{thrm:smooth:unconstrained:complexity}}\label{appendix:proof:thrm:smooth:unconstrained:complexity}

First, from \eqref{eq:proof:thrm:1:evn:1:1}, we have
\begin{align}\label{eq:proof:thrm:3:evn:1}
	F_{\lambda}(\bar{x}_{k+1})
	&\leq
		f(\hat{x}_k)
		+
		\frac{1}{2\lambda}
		\ltwo{ \bar{x}_k - \hat{x}_k - \stepsize g_k}^2		.
\end{align}
 We next follow \cite{DD19} and write $\hat{x}_k$ as 
\begin{align*}
	\hat{x}_k = \stepsize\lambda^{-1}\bar{x}_k - \stepsize\grad{f(\hat{x}_k)} + (1-\stepsize\lambda^{-1})\hat{x}_k,
\end{align*}
which is true since in this case $\grad{F_\lambda(\bar{x}_k)}=\grad{f(\hat{x}_k)}$.
Let $\delta=1-\stepsize\lambda^{-1}\in (0,1)$, we have
\begin{align*}
	\E[\ltwo{ \bar{x}_k - \hat{x}_k - \stepsize g_k}^2]
	&=
	\E[\ltwo{
		\delta(\bar{x}_k - \hat{x}_k) 
		-  
		\stepsize (g_k- \grad{f(x_k)})
		-
		\stepsize (\grad{f(x_k)} - \grad{f(\hat{x}_k)})
	}^2]
	\nonumber\\
	&=
	 \delta^2\E[\ltwo{ \bar{x}_k - \hat{x}_k}^2]
	 -
	 2\stepsize\delta\, \E\left[\InP{\bar{x}_k - \hat{x}_k}{\grad{f(x_k)} - \grad{f(\hat{x}_k)}}\right]
	\nonumber\\
	&\hspace{0.45cm}
	 +
	 \stepsize^2\E[\ltwo{\grad{f(x_k)} - \grad{f(\hat{x}_k)}}^2]
	 + \stepsize^2\sigma^2.
\end{align*}
Using  Young's inequality $|ab| \leq a^2/(2\theta)+\theta b^2/2$ for $a,b\in \R$ and $\theta>0$, we get
\begin{align*}
	 -2\stepsize\delta\, \E\left[\InP{\bar{x}_k - \hat{x}_k}{\grad{f(x_k)} - \grad{f(\hat{x}_k)}}\right]
	\leq
	\stepsize\delta\rho\, \E[\ltwo{ \bar{x}_k - \hat{x}_k}^2]
	+
	\stepsize\delta/\rho \,\E[\ltwo{\grad{f(x_k)} - \grad{f(\hat{x}_k)}}^2].
\end{align*}
Since $\grad{f}$ is $\rho$-Lipschitz, it follows that
\begin{align*}
	\E[\ltwo{ \bar{x}_k - \hat{x}_k - \stepsize g_k}^2]
	&\leq
		(\delta^2 + \rho\stepsize\delta)\E[\ltwo{ \bar{x}_k - \hat{x}_k}^2]
		+
		(\stepsize\delta/\rho+\stepsize^2)\E[\ltwo{\grad{f(x_k)} - \grad{f(\hat{x}_k)}}^2]
		+ 
		\stepsize^2\sigma^2
	\nonumber\\
	&\leq
		(\delta^2 + \rho\stepsize\delta)\E[\ltwo{ \bar{x}_k - \hat{x}_k}^2]
		+
		(\stepsize\delta\rho+\stepsize^2\rho^2)\E[\ltwo{x_k - \hat{x}_k}^2]
		+ 
		\stepsize^2\sigma^2	
	\nonumber\\
	&\leq
		(\delta^2 + 3\rho\stepsize\delta + 2\stepsize^2\rho^2)\E[\ltwo{ \bar{x}_k - \hat{x}_k}^2]
		+
		2(\stepsize\delta\rho+\stepsize^2\rho^2)\E[\ltwo{x_k - \bar{x}_k}^2]
		+ 
		\stepsize^2\sigma^2,			
\end{align*}
where we also used the inequality $\ltwo{a+b}^2\leq 2 \ltwo{a}^2+ 2 \ltwo{b}^2$ in the last step. We have
\begin{align*}
	\delta^2 + 3\rho\stepsize\delta + 2\stepsize^2\rho^2
	=
	1-\stepsize(2\lambda^{-1}-3\rho) + \stepsize^2\lambda^{-1}(\lambda^{-1}-3\rho+2\rho^2\lambda)
\end{align*}
Note that for $\lambda^{-1}\in[3\rho/2,2\rho]$, the term associated with $\stepsize^2$ is nonpositive, and hence
\begin{align}\label{eq:proof:thrm:3:evn:2}
	\E[\ltwo{ \bar{x}_k - \hat{x}_k - \stepsize g_k}^2]
	&\leq
		\E[\ltwo{ \bar{x}_k - \hat{x}_k}^2]
		-
		\stepsize(2\lambda^{-1}-3\rho)\E[\ltwo{ \bar{x}_k - \hat{x}_k}^2]		
	\nonumber\\
	&\hspace{0.45cm}
		+
		2(\stepsize\delta\rho+\stepsize^2\rho^2)\E[\ltwo{x_k - \bar{x}_k}^2]
		+ 
		\stepsize^2\sigma^2.				
\end{align}
Therefore, it follows from \eqref{eq:proof:thrm:3:evn:1} and \eqref{eq:proof:thrm:3:evn:2} that
\begin{align*}
	\E[F_{\lambda}(\bar{x}_{k+1})]
	\leq
		\E[F_{\lambda}(\bar{x}_{k})]
		-
		\stepsize\frac{2\lambda^{-1}-3\rho}{2\lambda^{-1}}\E[\ltwo{\grad{F_{\lambda}(\bar{x}_{k})}}^2]		
		+
		\frac{\stepsize\delta\rho+\stepsize^2\rho^2}{\lambda}\E[\ltwo{x_k - \bar{x}_k}^2]
		+ 
		\frac{\stepsize^2\sigma^2}{2\lambda}.
\end{align*}
We also have 
\begin{align*}
	\frac{\stepsize\delta\rho+\stepsize^2\rho^2}{\lambda}
	=
		\frac{\stepsize\rho}{\lambda} + \stepsize^2(\rho^2-\rho\lambda^{-1})
	\leq 
		\frac{\stepsize\rho}{\lambda}
	\leq 
		\frac{\stepsize}{\lambda^2},
\end{align*}
where the last two steps hold since $\lambda^{-1}\geq \rho$. Therefore, using the fact that $\xi=(1-\mmt)/\nu$, $\nu=\mmt/\stepsize$, $\bar{x}_k=x_k+\frac{1-\mmt}{\mmt}(x_k-x_{k-1})$, and $d_k=\frac{x_{k-1}-x_k}{\stepsize}$, we obtain
\begin{align}\label{eq:proof:thrm:3:evn:3}
	\E[F_{\lambda}(\bar{x}_{k+1})]
	\leq
		\E[F_{\lambda}(\bar{x}_{k})]
		-
		\stepsize\frac{2\lambda^{-1}-3\rho}{2\lambda^{-1}}\E[\ltwo{\grad{F_{\lambda}(\bar{x}_{k})}}^2]		
		+
		\frac{\stepsize\xi^2}{\lambda^2}\E[\ltwo{d_k}^2]
		+ 
		\frac{\stepsize^2\sigma^2}{2\lambda}.
\end{align}

Now, multiplying both sides of \eqref{eq:lem:smooth:W} by ${\xi^2}/\lambda^2$ and combining the result with \eqref{eq:proof:thrm:3:evn:3}, we obtain
\begin{align}\label{eq:proof:thrm:3:evn:4}
	\E\left[F_{\lambda}(\bar{x}_{k+1}) 
	+ 
	\frac{\xi^2}{\lambda^2}W_{k+1} 
	+ 
	\frac{\stepsize\xi^2}{\lambda^2}\ltwo{d_{k+1}}^2\right]
	&\leq
		\E\left[F_{\lambda}(\bar{x}_{k}) 
		+ 
		\frac{\xi^2}{\lambda^2}W_k 
		+ 
		\frac{\stepsize\xi^2}{\lambda^2}\ltwo{d_{k}}^2\right]		
	\nonumber\\
	&\hspace{-1.5cm}		
		-
		\stepsize\frac{2\lambda^{-1}-3\rho}{2\lambda^{-1}}\E[\ltwo{\grad{F_{\lambda}(\bar{x}_{k})}}^2]
		+
		\frac{(1+{8\nu\xi^2}/{\lambda})\stepsize^2\sigma^2}{2\lambda}.
\end{align}
Let $V_{k+1}$ denote the left-hand-side of \eqref{eq:proof:thrm:3:evn:4}, summing the result over $k=-1,\ldots,K-1$ gives
\begin{align*}
	 V_{K}
	\leq
		 {V_{-1}}
		-
		\frac{\stepsize_0}{\sqrt{K+1}}\frac{2\lambda^{-1}-3\rho}{2\lambda^{-1}}
		\sum_{k=0}^{K}
		\E[\ltwo{\grad{F_{\lambda}(\bar{x}_{k})}}^2]
		+
		\frac{(K+1)(1+{8\nu\xi^2}/{\lambda}) \sigma^2\stepsize^2}{2\lambda}.
\end{align*}
Note that 
\begin{align*}
	V_{-1} 
	= 
		F_{\lambda}(\bar{x}_{-1}) 
		+ 
		\frac{\xi^2}{\lambda^2}
		\left(2f(x_{-1}) + \frac{\varphi_{-1}}{\nu\stepsize^2} + \frac{\xi}{2}\ltwo{d_{-1}}^2\right)
		+
		\frac{\stepsize\xi^2}{\lambda^2}\ltwo{d_{-1}}^2
	\leq
		\left(1+\frac{2\xi^2}{\lambda^2}\right)f(x_0),
\end{align*}
where we used the facts that $x_{-1}=x_0$, $\varphi_{-1}=0$ (since $z_{-1}=0$), $d_{-1}=0$, and $F_{\lambda}(\bar{x}_{-1})=F_{\lambda}(\bar{x}_{0}) \leq f(x_0)$. 
Therefore, lower-bounding the left-hand-side by $(1 +2\xi^2/\lambda^2)f\opt$ and rearranging terms, we obtain
\begin{align}\label{eq:proof:thrm:3:evn:5}
	\E\left[\ltwo{\grad{F_{\lambda}(\bar{x}_{k^*})}}^2\right]
	\leq 
		\frac{2\lambda^{-1}}{2\lambda^{-1}-3\rho}\cdot
	\frac{
		(1 +2\xi^2/\lambda^2)\left(f(x_0)-f\opt\right)
		+ 
		\frac{(1+{8\nu\xi^2}/{\lambda}) \sigma^2\stepsize_0^2}{2\lambda}	
	}{
		\stepsize_0\sqrt{K+1}
	},	
\end{align}
where the last expectation is taken with respect to all random sequences generated by 
the method and the uniformly distributed random variable $k^*$. 
Since $\nu=1/\stepsize_0$ and $\xi=(1-\mmt)/\nu \leq 1/\nu$, hence 
\begin{align*}
	\E\left[\ltwo{\grad{F_{\lambda}(\bar{x}_{k^*})}}^2\right]
	\leq 
		\frac{2\lambda^{-1}}{2\lambda^{-1}-3\rho}\cdot
	\frac{
		(1 +2\stepsize_0^2/\lambda^2)\left(f(x_0)-f\opt\right)
		+ 
		\frac{(1+{8\stepsize_0}/{\lambda}) \sigma^2\stepsize_0^2}{2\lambda}	
	}{
		\stepsize_0\sqrt{K+1}
	}.	
\end{align*}
Finally, letting $\lambda=1/(2\rho)$ completes the proof.

\end{document}